\documentclass[12pt]{amsart}     

\usepackage{graphicx,morefloats}
\usepackage{epsfig,psfrag,caption}
\usepackage{graphicx,color,soul,verbatim,subfigure}
\usepackage{helvet,gbmacros,todonotes}
\usepackage[varg]{txfonts}

\theoremstyle{definition}

\theoremstyle{remark}

\newcommand{\goto}{\rightarrow}
\newcommand{\bigo}{{\mathcal O}}

\def\XXint#1#2#3{{\setbox0=\hbox{$#1{#2#3}{\int}$}
     \vcenter{\hbox{$#2#3$}}\kern-.5\wd0}}

\DeclareMathOperator{\res}{Res}
\DeclareMathOperator{\imag}{Im}

\DeclareMathOperator{\real}{Re}

\newcommand{\bea}{\begin{eqnarray}}
\newcommand{\eea}{\end{eqnarray}}
\newcommand{\bean}{\begin{eqnarray*}}
\newcommand{\eean}{\end{eqnarray*}}
\newcommand{\ba}{\begin{array}}
\newcommand{\ea}{\end{array}}
\newcommand{\beqs}{\begin{equation*}\begin{split}}

\newenvironment{choices}{\left\{ \begin{array}{ll}}{\end{array}\right.}
\newcommand\when{&\text{if~}}
\newcommand\otherwise{&\text{otherwise}}

\newenvironment{mat}{\left[\begin{array}{ccccccccccccccc}}{\end{array}\right]}
\newcommand\bcm{\begin{mat}}
\newcommand\ecm{\end{mat}}

\def\bigo{\mathcal O}
\def\eqref#1{(\ref{#1})}
\def\overl@ss#1#2{\vcenter{\offinterlineskip
        \ialign{$\m@th#1\hfil##\hfil$\crcr#2\crcr<\crcr } }}
\def\overgr@at#1#2{\vcenter{\offinterlineskip
        \ialign{$\m@th#1\hfil##\hfil$\crcr#2\crcr>\crcr } }}
\def\gl{\mathrel{\mathpalette\overl@ss>}}
\def\lg{\mathrel{\mathpalette\overgr@at<}}
\def\d{\mathrm{d}}

\def\Real{\mathbb{R}}
\def\Complex{\mathbb{C}}
\def\Re{\mathop{\rm Re}\nolimits}
\def\Im{\mathop{\rm Im}\nolimits}
\def\arg{\mathop{\rm arg}\nolimits}

\def\pvint{\int\kern-0.94em-\kern0.2em}
\let\^=\hat
\let\==\bar
\let\@=\mathbf

\let\le=\leqslant
\let\ge=\geqslant
\def\d{\mathrm{d}}
\def\e{\mathrm{e}}

\def\~#1{\tilde{#1}}
\makeatother

\def\be{\begin{equation}}
\def\ee{\end{equation}}
\def\bse{\begin{subequations}}
\def\ese{\end{subequations}}

\def\pvint{\mathop{\int\kern-0.84em-\kern0.2em}\limits}
\def\erf{\mathop{\rm erf}\nolimits}

\def\res{\mathrm{res}}
\def\heat{\mathrm{heat}}
\def\schr{\mathrm{schr}}

\def\reg{\mathrm{reg}}
\numberwithin{figure}{section}

\let\eref=\eqref
\newenvironment{proof-univ}{\paragraph{Proof of Theorem~\ref{t:universality}:}}{\hfill$\square$}

\begin{document}                        


\title{Gibbs phenomenon for dispersive PDEs on the line}

\author{Gino Biondini}
\address{Department of Mathematics,
State University of New York at Buffalo,
Buffalo, NY 14260}
\author{Thomas Trogdon}
\address{Courant Institute of Mathematial Sciences,
New York University,
251 Mercer St.,
New York, NY, USA}

\begin{abstract}
We investigate the Cauchy problem for linear, constant-coefficient evolution PDEs on the real line with discontinuous initial conditions (ICs) in the small-time limit.  The small-time behavior of the solution near discontinuities is expressed in terms of universal, computable special functions.  We show that the leading-order behavior of the solution of dispersive PDEs near a discontinuity of the ICs
is characterized by Gibbs-type oscillations and gives exactly the Wilbraham-Gibbs constant.
\end{abstract}

\maketitle   


\setcounter{tocdepth}{1} 
\tableofcontents


\section{Introduction}
\label{s:intro}

\begin{figure}[t]
\centerline{\includegraphics[width=0.9\linewidth]{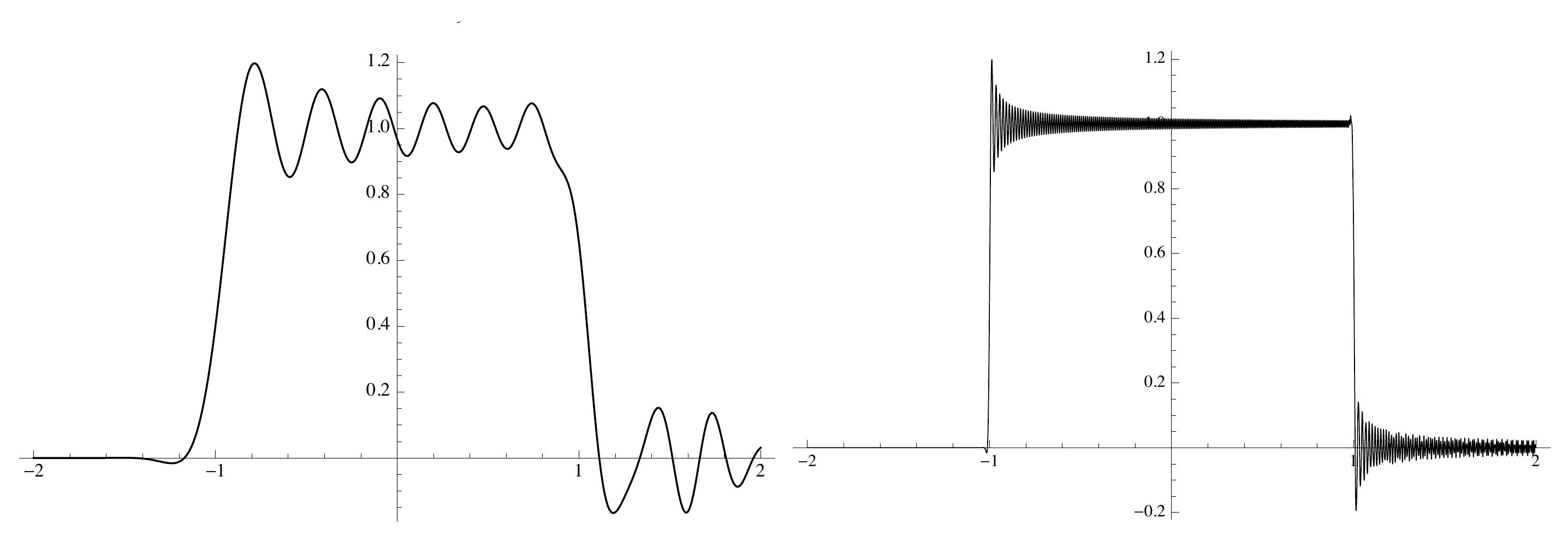}}
\caption{The solution of \eqref{e:pde} with $\omega(k) = k^5$ and IC $q_o(x) = 1$ if $|x| \leq 1$ and $q_o(x) = 0$ otherwise.  
Left: $t=10^{-6}$. Right: $t= 10^{-12}$.
The solution exhibits the Gibbs phenomenon, as discussed in detail in Section~\ref{s:generaldispersion}.}
\label{f:fifth} 
\end{figure}

The Gibbs phenomenon is the well-known behavior of the Fourier series of a piecewise continuously differentiable periodic function at a jump discontinuity.
Namely, the partial sums of the Fourier series have large oscillations near the jump, 
which typically increase the maximum of the sum above that of the function itself
\cite{carslaw,gibbs}. 
Moreover, the overshoot does not subside as the frequency increases, but instead approaches a finite limit.
The Gibbs phenomenon is typically viewed as a numerical artifact in the numerical representation of a function due to truncation.  
Here we view the Gibbs phenomenon as a product of \emph{non-uniform convergence}.
Namely, the partial sums of the Fourier series are analytic and converge to a discontinuous function, 
and hence, this convergence must be non-uniform in any neighborhood of the discontinuity.
In turn, this gives rise to highly oscillatory behavior.
Keeping this view of the Gibbs phenomenon in mind, we show in this work that 
\textit{the solution of dispersive PDEs with discontinuous initial conditions (ICs) exhibit the Gibbs phenomenon for short times}.
E.g., Fig.~\ref{f:fifth} shows a solution of \eqref{e:pde} with $\omega(k) = k^5$ for short times. 

Specifically, we consider initial value problems (IVPs) of the form
\begin{gather}
i q_t - \omega(-i\partial_x) q = 0,
\label{e:pde}
\\
\noalign{\noindent with $\omega(k) = \omega_n k^n + \bigo(k^{n-1})$ and with IC}
q(x,0) = q_o(x).
\label{e:IC}
\end{gather}
Unless otherwise stated, we assume the $\omega : \mathbb R \goto \mathbb R$.
It is well known that, for hyperbolic PDEs, 
the discontinuities of the IC travel along characteristics~ 
\cite{couranthilbert,evans}.  
For dispersive and diffusive PDEs, in contrast, 
even if the ICs are discontinuous, the solution of the IVP is typically classical $\forall x\in\Real$ as long as $t>0$ and the IC has sufficient decay as $|x|\to\infty$.
But an obvious question is: \emph{What does the solution actually look like as $t \downarrow 0$?}
Answering this question is useful for many reasons. 
For example:
(i) to evaluate asymptotics for linear and nonlinear problems \cite{Taylor2006},
(ii) to build/test numerical integrators, and
(iii) to understand the behavior of initial boundary value problems (IBVPs).  Surprisingly, however, 
while the smoothing effects of diffusion are well known, this perspective on dispersive regularization 
is not as well characterized in the literature to the best of our knowledge.  One of the central messages of this work is that: 
\textit{The slow decay of the Fourier transform of $q_o$ as $k\to\infty$ affects the short-time asymptotics of the solution $q(x,t)$}.

Let us briefly elaborate on item~(iii) above.
One of the original motivations for this work was the study of corner singularities 
in IBVPs \cite{boydflyer,flyerfornberg2003a,flyerfornberg2003b,flyerswarztrauber}.  
The issue at hand is the following.  
Consider \eqref{e:pde}, with $n=2$, posed on the domain $D = (0,\infty) \times (0,T)$ so that one has to also specify boundary data at $x = 0$, say $q(0,t) = g_0(t)$.  
The smoothness of $q(x,t)$ in $\overline D$ is restricted not only by the smoothness (and decay) of $q_o(x)$ and $g_0(t)$ but by the compatibility of these two functions at $x = t = 0$, \emph{i.e.}, to first-order,
$q_o(0) = g_0(0)$.  
(Higher-order conditions are found by enforcing the PDE holds at the corner of the domain.)
When compatibility fails at some order, a corner singularity is present.  
One would obviously like to characterize the effect of such a corner singularity on the solution of an IBVP.
It soon became clear, however, that in order to do so, one needs to fully understand the behavior of IVPs with discontinuous ICs.

The outline of this work is the following.  
In Section~\ref{s:summary} we summarize our fundamental results concerning both the smoothness of solutions and their short-time behavior. 
In Section~\ref{s:onediscontinuity} we perform the asymptotic analysis in the case of a single discontinuity in the IC $q_o$.  
There we identify the special functions that describe the Gibbs-like  behavior.  
Such functions are generalizations of the classical special functions, and are computable with similar numerical methods.
In Section~\ref{s:generaldispersion} we display some sample solutions, we discuss their Gibbs-like behavior,
we further study the properties of the special functions and we establish a precise connection with the classical Gibbs phenomenon.
In Section~\ref{s:discontinuousderivatives} we treat the case where $q_o'$ has one jump discontinuity.  
In Section~\ref{s:multiplediscontinuities} we present our general result, which allows for multiple discontinuities in $q_o$ itself or in any of its derivatives.  
A full asymptotic expansion is derived near, and away from, the singular (i.e., non-smooth) points of $q_o$.
Section~\ref{s:specialfunctions} contains additional details on the analysis and numerical computation of the special functions considered here.
Finally, Section~\ref{s:discussion} concludes this work with a discussion of the results and some final remarks.

Further details and technical results are relegated to four appendices.  
In Appendix~\ref{a:wellposedness} we review some well-known results about well-posedness of the IVP.
In Appendix~\ref{a:regularity} we prove the result (stated in Section~\ref{s:regularity}) 
concerning the classical smoothness of the solution for $t > 0$, using the method of steepest descent for integrals. 
Appendix~\ref{a:residual} contains technical results for determining the order of the error terms in our short-time expansions.  
Finally, in Appendix~\ref{a:approx} we study the robustness of the Gibbs phenomenon by analyzing
the behavior of solutions whose ICs are a small perturbation of a discontinuous function.

The function $\omega(k)$ is referred to as the \textit{dispersion relation} of the PDE \eref{e:pde}.
Throughout this work, we will take the dispersion relation $\omega(k)$ to be polynomial.
Note that one can always remove constant and linear terms from $\omega(k)$ by performing a phase rotation and a Galilean transformation, respectively.
That is, without loss of generality one can take
\[
\omega(k) = \sum_{j=2}^n \omega_j k^j\,. 
\label{e:dispersion}
\]
We will assume that this has been done throughout this work.

\section{Summary of results}
\label{s:summary}

This section contains a brief summary of our main results.  Our summarized results concern regularity, the Gibbs phenomenon and asymptotics for our special functions.  Another one of our main results (Theorem~\ref{Thm:Main}) is not summarized here due to its complexity: It gives the full expansion of the solution  of \eqref{e:pde} for short times.

\subsection{Regularity results for linear evolution PDEs}
\label{s:regularity}

We begin this section by referring to Appendix~\ref{a:wellposedness} for the required definitions and classical results concerning the well-posedness of \eqref{e:pde} for $q_o \in L^2(\mathbb R)$, where
\begin{align}\label{e:soln-form}
  q(x,t) &= \frac{1}{2 \pi} \int_{-\infty}^\infty e^{i\theta(x,t,k)} \hat q_o(k) \d k,\\
  \theta(x,t,k) &= kx - \omega(k) t. \label{e:thetadef}
\end{align}
Two properties can be readily seen:
\begin{enumerate}
\item[1.]
$q(\cdot,t) \rightarrow q(\cdot,0)$ in $L^2(\mathbb R)$ as $t \downarrow 0$, 
\item[2.] 
if $q_o \in H^1(\mathbb R)$ then $q(\cdot,t) \rightarrow q(\cdot,0)$ uniformly as $t \downarrow 0$.
\end{enumerate}

More delicate questions can be asked about pointwise behavior in the short-time limit, however.
In particular, 
Sj\"olin \cite{Sjolin} showed that when $\omega(k) = k^{2m}$, $m = 1,2 \ldots$, and $q_o \in H^s(\mathbb R)$ with compact support for $s \geq 1/4$, $\lim_{t \downarrow 0} q(x,t) = q(x,0)$ for a.e. $x \in \mathbb R$.
This result was generalized in \cite{kenigponcevega} for general $\omega(k)$ without the assumption of compact support.
(See also \cite{Lee,Tao,Vega}.)
The results that follow will only demonstrate a.e.\ convergence for a subset of $H^{1/4}(\mathbb R)$.
On the other hand, the short-time expansion that we will provide in the following sections is new.

Interesting questions related to the regularity of the solution can also be asked.
As is noted in \cite{TaoCBMSbook}, when $\omega(k) = k^2,~k^3$ the solutions are easily seen to be continuous for $t > 1$ provided $q \in L^2 \cap L^1(\mathbb R)$.  Furthermore the $L^\infty(\mathbb R)$ norm of $q(\cdot,t)$ decays in time.
A Strichartz-type result was provided in \cite{kenigponcevega} showing, in particular, the space-time estimate $\|q\|_{L^8(\mathbb R^2)} \leq C \|q_o\|_{L^2(\mathbb R)}$ for $\omega(k) =k^3$.  In Appendix~\ref{s:regularity}
we prove results of a classical nature concerning the regularity of the solution:

\begin{theorem}[Regularity]
\label{t:differentiable}
Let $\omega(k)$ be as in~\eref{e:dispersion} and $q(x,t)$ as in~\eref{e:soln-form},
with $q_o \in L^2(\mathbb R) \cap L^1(\mathbb R, (1+|x|)^\ell \d x)$.
\begin{itemize}
\item[(i)] If
\begin{align*}
\ell \geq \frac{2m -n + 2}{2(n-1)},
\end{align*}
$q(x,t)$ is differentiable $m$ times with respect to $x$ for $t >0$ and $\partial_x^m q(x,t)$ is continuous as a function of $x$ and $t$.
\item[(ii)] If 
\begin{align*}
\ell \geq \frac{2jn -n +2}{2(n-1)},
\end{align*}
$q(x,t)$ is differentiable $j$ times with respect to $t$ for $t > 0$ and $\partial_t^j q(x,t)$ is continuous as a function of $x$ and $t$.
\end{itemize}
\end{theorem}

\begin{corollary}[Classical solution]
\label{c:classical}
Under the same hypotheses as in Theorem~\ref{t:differentiable}, 
if
\begin{align*}
\ell \geq \mathfrak C_n \triangleq \frac{n+2}{2(n-1)},
\end{align*}
the $L^2$ solution of the IVP is classical (see Definition~\ref{d:classical}) for $t > 0$.
\end{corollary}

The importance of these results from the perspective of this paper is that if we can guarantee that the solution is smooth for $t > 0$ and if the IC is not smooth then we can guarantee that the limit $t \downarrow 0$ is a singular one: It forces the breakdown of smoothness.  The last regularity result concerns the integrability of solutions.

\begin{corollary}[Loss of integrability]
\label{c:non-integrable}
Let $\omega(k)$ be as in~\eref{e:dispersion} and $q(x,t)$ as in~\eref{e:soln-form},
with $q_o \in L^1 \cap L^2(\mathbb R)$.  Assume $q_o$ has at least one jump discontinuity\footnote{To be precise, we assume $\real q_o(x) = \limsup_{\delta \downarrow 0} \int_{|y-x| < \delta} \real q_o(y) \d y$ and $\imag q_o(x) = \limsup_{\delta \downarrow 0} \int_{|y-x| < \delta} \imag q_o(y) \d y$.}.  Then $q(\cdot,t) \not\in L^1(\mathbb R)$ for any $t > 0$.
\end{corollary}
\begin{proof}
  Assume $\tilde q_o \triangleq q(\cdot, t) \in L^1(\mathbb R)$ for some $t > 0$.  Then take this as an initial condition for the PDE with $\omega(k)$ replaced with $-\omega(k)$ and find its solution $\tilde q(x,t)$.  Then $\tilde q(x,t)$ should be continuous as function of $x$ by Theorem~\ref{t:differentiable} but this is a contradiction as uniqueness ensures $q_o = \tilde q(\cdot, t)$ and $q_o$ is discontinuous. 
\end{proof}

\subsection{Short-time behavior}

To explain two aspects of the short-time behavior we state some theorems.  Define
\[
I_{\omega,0}(y,t) \triangleq \frac1{2\pi}\,\int_C \e^{iky- i\omega(k)t}\frac{\d k}{ik}\,,
\]
where $C$ is a contour in the closed upper-half plane that runs along the real axis but avoids $k = 0$.

\begin{theorem}[Leading-order universality]\label{t:universality}
  Assume $q_o \in L^2(\mathbb R)$ and there exists $c_0 = -\infty < c_1 < c_2 < \cdots < c_N < c_{N+1} = \infty$ such that the restriction $q|_{(c_i,c_{i+1})}$ has one derivative in $L^2((c_i,c_{i+1}))$ for each $i = 0,1,\ldots, N$.  Then if $[q_o(c_i)] \triangleq q_o(c_i^+) - q_o(c_i^-) \neq 0$, there exists a constant $q_{c_i}$ such that
  \begin{align*}
    \lim_{t \downarrow 0} \frac{ q(c_i + x |\omega_nt|^{1/n},t) - q_{c_i} }{[q_o(c_i)]} = I_{\omega_n,0}(x,1), \quad \omega_n(k) = e^{i \arg(\omega_n)} k^n,
  \end{align*}
  uniformly for $x$ in a bounded set.
\end{theorem}
This is interpreted as a universality theorem because, after proper rescaling, the solution is the same independent of both the initial condition and the lower terms in the dispersion relation.  It is proved in Section~\ref{s:onediscontinuity}.  Because of the differential equation, \eqref{e:diff-eq} satisfied by $I_{\omega_n,0}(x,1)$ we have that
\textit{ the leading-order behavior\footnote{One can generalize this with appropriate scaling when any derivative of $q_o$ is discontinuous but we do no pursue this further here.  This gives a universality statement involving $I_{\omega_n,m}$.} of the solution near a discontinuity is governed by
  a similarity solution expressed in terms of classical special functions.}

The non-uniform convergence of $q(x,t)$ to $q_o(x)$ as $t \downarrow 0$ when $q_o(x)$ is discontinuous at $x = c$ generically results in a so-called overshoot value --- the amount by which $q(x,t)$ over (or under) approximates $q_o(c^\pm)$, see Figure~\ref{f:fifth}. We relate the behavior of the overshoot near this region of non-uniform convergence as $t \downarrow 0$ to the Gibbs phenomenon with the following theorems.  The first is a restatement of the results of Wilbraham and Gibbs (\cite{wilbraham} and \cite{gibbs}):
\begin{theorem}[Gibbs phenomenon]
Consider the Fourier series approximation of
\begin{align*}
f(x) = \begin{choices}1, \when |x| \leq 1,\\
0, \otherwise,
\end{choices}
~~\text{given by}~~~
S_n[f](x) =\sum_{k=-n}^{n} \frac{4 \sin \frac{k \pi}{2}}{k \pi} \mathrm{e}^{\frac{\mathrm {i} k x \pi}{2}}.
\end{align*}
Then for any $\delta > 0$
\begin{align*}
\lim_{n \goto \infty} \sup_{|x\pm 1| \leq \delta} S_n[f](x) &= 1 + \mathfrak g,\\
\lim_{n \goto \infty} \inf_{|x\pm 1| \leq \delta} S_n[f](x) &= -\mathfrak g,
\end{align*}
where
\begin{align*}
\mathfrak g = \frac{1}{\pi} \int_0^\pi \frac{\sin z}{z} \,\d z  - \half \approx 0.089490\ldots.
\end{align*}
\end{theorem}
In this context our results give:
\begin{theorem}[Gibbs phenomenon on the line]\label{t:gibbs}
Let $q_n(x,t)$ be the solution of $i q_t - (- i \partial_x)^n q = 0$ with
\begin{align*}
q(x,0) = \begin{choices}1, \when |x| \leq 1,\\
0, \otherwise.
\end{choices}
\end{align*}
Then for any $\delta > 0$
\begin{align*}
\lim_{n \goto \infty} \lim_{t \downarrow 0} \sup_{|x\pm 1| \leq \delta} \real q_n(x,t) &= 1 + \mathfrak g,\\
\lim_{n \goto \infty} \lim_{t \downarrow 0} \inf_{|x\pm 1| \leq \delta} \real q_n(x,t) &= -\mathfrak g,\\
\lim_{n \goto \infty} \lim_{t \downarrow 0} \sup_{|x\pm 1| \leq \delta} \imag q_n(x,t) &= 0,\\
\lim_{n \goto \infty} \lim_{t \downarrow 0} \inf_{|x\pm 1| \leq \delta} \imag q_n(x,t) &= 0.
\end{align*}
\end{theorem}

  One does not have to take $\omega(k) = k^n$ in the previous theorem: It follows for general $\omega(k)$ provided the coefficients are appropriately controlled.  One such example is
  \begin{align*}
    \omega(k) = k^n + \sum_{j=n-m}^{n-1} c_{j,n} k^j,
  \end{align*}
  where $-C \leq c_{j,n} \leq C$ are real and $m$ is fixed.  Furthermore, there is an analog of this theorem that holds for general data as in Theorem~\ref{t:universality}.  This phenomenon is explored in greater depth in Section~\ref{s:gibbs}.

\textit{The Gibbs-like oscillations represent the real behavior of the solution of dispersive PDEs, and are not a numerical artifact.}
In other words, Fig.~\ref{f:schrodinger} (as well as Fig.~\ref{f:stokes} and the figures in Section~\ref{s:specialfunctions}) 
are not a result of truncation error!  This fact has important consequences for the numerical solution of dispersive PDEs, particularly, in finite-volume methods where a so-called Riemann problem must be solved.

\subsection{Asymptotics of $I_{\omega,m}$}

The previous results rely on the asymptotic analysis of the function $I_{\omega,m}(x,t)$ as $t \downarrow 0$ or as $|x| \goto \infty$ for fixed $t > 0$. We also define the kernel $K_t(x)$ by
\begin{align*}
  q(x,t) = \frac{1}{2 \pi} \int_{-\infty}^\infty e^{ikx - i \omega(k)t} \hat q_o(k) \d k = \int_{-\infty}^\infty K_t(x-y) q_o(y) \d y.
\end{align*}
In Appendix~\ref{a:regularity} we use the method of steepest descent for integrals to derive precise asymptotics of $I_{\omega,m}$ and $K_t(x) = I_{\omega,-1}(x,t)$.  First, we rescale the integral
\begin{align}
\begin{split}
  I_{\omega,m}(x,t) &= \frac{1}{2 \pi}\sigma^m \left(\frac{|x|}{t}\right)^{-m/(n-1)} \int_{C} e^{X(iz-i\omega_n \sigma^n z^n-i R_{|x|/t}(z))} \frac{\d z}{(iz)^{m+1}},\\
  \sigma &= \sign(x), \quad k = \sigma (|x|/t)^{1/(n-1)}z,\\
R_{|x|/t}(z) &\triangleq \sum_{j=2}^{n-1} \omega_j \left(\frac{|x|}{t}\right)^{\frac{j-n}{n-1}} (\sigma z)^j, \quad X \triangleq |x| \left( \frac{|x|}{t} \right)^{1/(n-1)},\\
\Phi_{|x|/t}(z) &= iz-i\omega_n \sigma^n z^n-i R_{|x|/t}(z),
\end{split}
\end{align}
Then define $\{z_j\}_{j=1}^{N(n)}$ to be the the solutions of $\Phi'_{|x|/t}(z) = 0$ in the closed upper-half plane.  Finally, define $\theta_j$ to be the direction at which the path of steepest descent leaves $z_j$ with increasing real part.  

\begin{theorem}\label{Thm:Kernel}
As $|x/t| \rightarrow \infty$
\begin{align*}
  I_{\omega,m}(x,t) &= -i \Res_{k=0}  \left(\frac{ e^{ikx-i\omega(k)t}}{(ik)^{m+1}}\right) \chi_{(-\infty,0)}(x) \\
  &+  \frac{\sigma^m |x|^{-1/2}}{\sqrt{2 \pi}}\left(\frac{|x|}{t}\right)^{-\frac{m+1/2}{n-1}} \sum_{j=1}^{N(n)}   \frac{e^{X \Phi_{|x|/t}(z_j) +i\theta_j}}{(i z_j)^{m+1}} \frac{1}{|\Phi_{|x|/t}''(z_j)|^{1/2}}\left(1 + \bigo\left(|x|^{-1} \left( \frac{|x|}{t} \right)^{-1/(n-1)}\right)\right).
\end{align*}
Hence:
\begin{itemize}
\item  For fixed $t> 0$ as  $|x| \rightarrow \infty$
\begin{align}\label{e:largespace}
K^{(m)}_t(x) \leq c \left\{ \begin{array}{lr} |x|^{\frac{2m-n+2}{2(n-1)}}, & n \text{ is even},\\
\\
|x|^{\frac{2m-n+2}{2(n-1)}}, & n \text{ is odd}, ~~ \omega_n x > 0,\\
\\
|x|^{-M} \text{ for all } M >0, & n \text{ is odd}, ~~ \omega_n x < 0,
\end{array}\right.
\end{align}
where $c$ depends on $m$, $t$ and $n$.
\item  For $|x| \geq \delta > 0$ and $m \geq 0$ as $t \rightarrow 0^+$
\begin{align}\label{e:smalltime}
I_{\omega,m}(x,t) &= -i \Res_{k=0}  \left(\frac{ e^{ikx-i\omega(k)t}}{(ik)^{m+1}}\right) \chi_{(-\infty,0)}(x) + \bigo\left(t^{\frac{m+1/2}{n-1}} |x|^{-\frac{2m+2n}{2(n-1)}} \right).
\end{align}

\end{itemize}
\end{theorem}

\section{Short-time asymptotics: discontinuous ICs}\label{s:onediscontinuity}

Recall that the above representation for the weak solution \eref{e:soln-form} of the IVP is valid as long as the IC $q_o(x)$ belongs to $L^2(\Real)$.  
We first consider initial data with a single discontinuity.
For now we will assume that $q_o$ satisfies the following properties:

\begin{assume}
Let
\begin{itemize}\label{Assume:One}
\item$q_o \in L^2(\Real)$,
\item$[q_o(c)] \triangleq q_o(c^+)-q_o(c^-) \neq 0$,
\item$q'_o$ exists on $(-\infty,c)\cup (c,\infty)$,
\item$q'_o \in L^q(-\infty,c) \cap L^q(c,\infty)$ for some $1 < q < \infty$, and
\item $q_o$ is compactly supported.
\end{itemize}
\end{assume}

In later sections we will discuss the effect of discontinuities in the derivatives of the IC and we will remove the condition of compact support.
The phenomenon we wish to investigate here is the following.  
The solution is classical for $t > 0$, but converges to a discontinuous function as $t \rightarrow 0$.  
Thus, the limit generally exists in $L^2(\mathbb R)$ but must fail to be uniform.  

To derive  an expansion for the solution for short times it is convenient to integrate the definition~\eref{e:ft-pair} of the Fourier transform by parts:
\begin{gather}
\^q_o(k) = \bigg(\!\! \int_{-\infty}^c + \int_c^\infty \bigg) \e^{-ikx}q_o(x)\,\d x
 = \frac1{ik}\e^{-ikc}[q_o(c)] + \frac1{ik} F(k)\,,
\label{e:q0hatparts} 
\\
\noalign{\noindent where}
F(k) = \bigg(\!\! \int_{-\infty}^c + \int_c^\infty \bigg) \e^{-ikx}q_o'(x)\,\d x\,, \quad [q_o(c)] = q_o(c^+) - q_o(c^-).
\label{e:Fdef}
\end{gather}
In Appendix~\ref{a:residual} we discuss the properties of $F(k)$.
Note that both terms in the right-hand side (RHS) of~\eref{e:q0hatparts} are singular at $k=0$, but their sum $\^q_o(k)$ is not.  Inserting~\eref{e:q0hatparts} in the reconstruction formula~\eref{e:ft-pair} for the solution of the IVP yields:
\[
q(x,t) =  \frac1{2\pi}[q_o(c)]\pvint_\Real \e^{i(k(x-c)-\omega(k)t)}\frac{\d k}{ik} + \frac1{2\pi}\pvint_\Real \e^{i\theta(x,t,k)} F(k)\frac{\d k}{ik}\,.
\label{e:qdecomp0}
\]
where $\pvint$ denotes the principal value (p.v.) integral.
The principal value sign is now needed because each of the integrands in~\eref{e:qdecomp0} is separately singular
at $k=0$.
Of course, one could have chosen other ways to regularize the singularity, 
and the final result for $q(x,t)$ is independent of this choice.

\begin{figure}[t!]
\kern+\smallskipamount
\centerline{\includegraphics[width=0.5\textwidth]{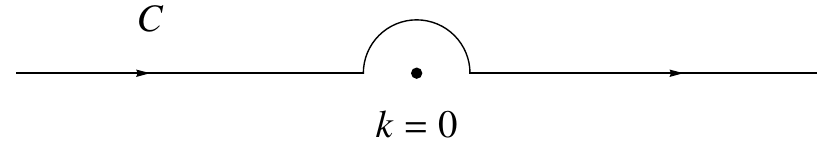}}
\kern-\smallskipamount
\caption{The integration contour $C$ for the evaluation of the principal value integral in \protect\eref{e:pvintcontour}. We assume the radius of the semi-circle is less than 1.}
\label{f:contour}
\end{figure}

We next show that the second term in the RHS of~\eref{e:qdecomp0} is continuous as a function of $x$ for all $t\ge0$, 
while the first term yields the dominant behavior in the neighborhood of the discontinuity at short times.
More precisely, we can write the p.v.\ integral in~\eref{e:qdecomp0} as:
\[
\pvint_\Real \!f(k)\,\d k = \int_C f(k)\,\d k + \pi i\Res_{k=0}[f(k)]
\label{e:pvintcontour}
\]
where $C$ is the contour shown in Fig.~\ref{f:contour}.
Recall
\[
I_{\omega,0}(y,t) \triangleq \frac1{2\pi}\,\int_C \e^{i[ky-\omega(k)t]}\frac{\d k}{ik}\,.
\label{e:I0def}
\]
(The reason for the subscript ``0'' will become apparent later on when we generalize these results to discontinuities in the higher derivatives.)
Also, define
\begin{align*}
q_c &= \half[q_o(c)] + \frac1{2\pi}\pvint_\Real \e^{ikc}F(k)\frac{\d k}{ik}, \quad  q_\res(y,t) = \frac1{2\pi}\int_\Real \e^{ikc}\frac{\e^{i\theta(y,t,k)}-1}{ik}\,F(k)\,\d k\,.
\end{align*}
Recalling $\Res\nolimits_{k=0}(\e^{i\theta(x{-}c,t,k)}/k)=1$,
we then write the decomposition~\eref{e:qdecomp0} as
\begin{equation}
q(x,t) = q_c + [q_o(c)] I_{\omega,0}(x-c,t) + q_\res(x-c,t)\,. 
\label{e:qdecomp}
\end{equation}
Note that the principal value is not needed on $q_\res(y,t)$, because the integrand is continuous at $k=0$.
Note also that the above decomposition holds for an arbitrary dispersion relation $\omega(k)$.

\begin{figure}[t!]
\centerline{\includegraphics[width=0.475\textwidth]{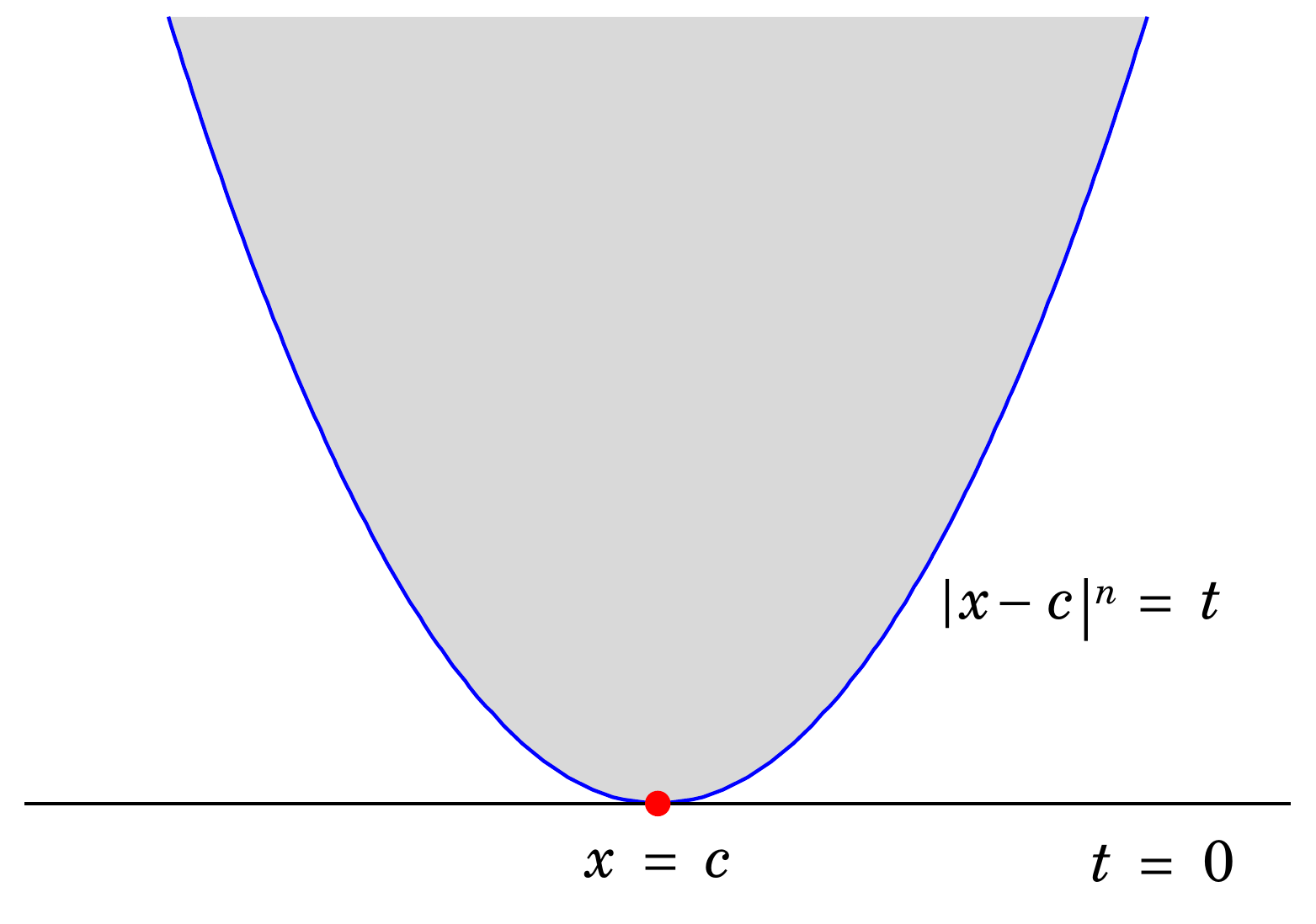}}
\caption{The regularization region (in gray) around a discontinuity in the IC.}
\label{f:regularizationregion}
\end{figure}

Importantly, each of the three terms in~\eref{e:qdecomp} are individually a solution of the PDE~\eref{e:pde}.  
However, each of them provides a different type of contribution.  
Indeed, a closer look allows the following interpretation of these pieces:
\begin{enumerate}
\item[(i)] $q_c$ represents a constant offset.
\item[(ii)] $[q_o(c)] I_{\omega,0}(y,t)$ characterizes the dominant behavior near the jump discontinuity.
The detailed properties of $I_{\omega,0}(y,t)$ are discussed in Appendix~\ref{a:regularity}.
In particular, Theorem~\ref{Thm:Kernel} implies
\[
\lim_{y\to\infty} I_{\omega,0}(y,t) = 0\,,\qquad
\lim_{y\to-\infty}I_{\omega,0}(y,t) = -1\,.
\label{e:I0limits}
\]
Note also that $\lim_{t \downarrow 0}I_{\omega,0}(0,t)\neq 0$. 
\item[(iii)]$q_\res(c,0)=0$, and $q_\res(x,t)$ is H\"older continuous and vanishes at $(x,t)=(c,0)$ for $t \geq 0$.
\end{enumerate}
One can look at the last item essentially as a trivial consequence of the first two, because 
the offset value and the jump behavior are all captured by the first and second contribution, respectively.
In practice, however, the proof is done in the reverse.   
Namely, in Appendix~\ref{a:residual} we prove~(iii) and we obtain precise estimates for the behavior of 
$q_\res(x-c,t)$ near $(x,t)=(c,0)$.  
More precisely, we show that, for $\|F\|_{L^p(\mathbb R)}< \infty$,
\begin{align}
\label{e:on-estimate}
q_\res(x-c,t) = \bigo( |x-c|^{1/p} + |t|^{1/(np)}).
\end{align}
The error term in the above short-time expansion is consistent as $t\to0$ as long as $|x-c|^n = O(t)$.
That is, the above expansion is valid in the region $|x-c|^n \le Ct$ (for some $C>0$)
in the neighborhood of a discontinuity $c$. 
We call such region the \textit{regularization region}. 
Such a region is illustrated in Figure~\ref{f:regularizationregion}.

One may also wish to understand the behavior of the solution in the short-time limit away from the singularity. 
Of course, to leading order, we expect it to be unaffected by the singularity and to limit pointwise to the IC.  
To prove that this is indeed the case, one must derive an estimate for the error term.  
The asymptotics of $I_{\omega,m}(x-c,t)$ can be fully characterized, see Theorem~\ref{Thm:Kernel}.  
The relevant behavior for the present purposes is
\begin{align*}
I_{\omega,0}(x-c,t) &= -\chi_{(-\infty,0)}(x-c) +  \bigo(t^{1/(2(n-1))})
\end{align*}
as $t\to0$ uniformly in the region $|x-c| \geq \delta > 0$.
Here and below, $\chi_R(y)$ is the characteristic function of a set $R$. (Namely, $\chi_R(y)=1$ for $y\in R$ and 
$\chi_R(y)=0$ otherwise.)
We then have 
\begin{align*}
q(x,t) = [q_o(c)](\half-\chi_{(-\infty,0)}(x-c)) + \frac{1}{2\pi} \pvint_{\mathbb R} e^{i\theta(x,t,k)} F(k) \frac{\d k}{ik}
+ \bigo(t^{1/(2(n-1))}).
\end{align*}
The relevant tool for the characterizing the limiting behavior of the rest of the solution is Lemma~\ref{Lemma:ResEst-HO}. 
From that result, \eqref{e:q0hatparts} and the above discussion it follows that
\begin{equation*}
q_0(x) 
= [q_o(c)](\half-\chi_{(-\infty,0)}(x-c)) +\frac{1}{2 \pi} \pvint_{\mathbb R} e^{ikx} F(k) \frac{\d k}{ik}.
\end{equation*}
Therefore, for $|s-c| \geq \delta > 0$ and $\|F\|_{L^p(\mathbb R)}< \infty$, we have
\begin{align}
\label{e:away-estimate}
q(x,t) = q_0(s) + \bigo( |x-s|^{1/p} + |t|^{1/(np)}+ |t|^{1/(2(n-1))}).
\end{align}
These observations also allow us to prove Theorem~\ref{t:universality}.

\begin{proof}[Proof of Theorem~\ref{t:universality}]
  Under Assumption~\ref{Assume:One}
  \begin{align*}
    \lim_{t \downarrow 0} \frac{ q(c + x |\omega_n|^{1/n}t^{1/n},t) - q_c }{[q_o(c)]} = \lim_{t \downarrow 0} I_{\omega_n,0}(x|\omega_n|^{1/n}t^{1/n},t),
  \end{align*}
  follows directly from \eqref{e:on-estimate}.  Then
  \begin{align*}
    I_{\omega_n,0}(x|\omega_n|^{1/n}t^{1/n},t) = \frac{1}{2\pi} \int_C e^{i(k |\omega_n|^{1/n}t^{1/n}) x - i \arg(\omega_n) ( |\omega_n|^{1/n} k t^{1/n})^n - r(k)t} \frac{\d k}{ik},
  \end{align*}
  where $r(k)$ is a polynomial of degree at most $n-1$.  Using $k t^{1/n}|\omega_n|^{1/n} \mapsto k$, and redeforming $C$, we have
  \begin{align*}
    I_{\omega_n,0}(x|\omega_n|^{1/n}t^{1/n},t) = \frac{1}{2\pi} \int_C e^{i k x - i \arg(\omega_n) k^n - r(k|\omega_n|^{-1/n}t^{-1/n})t} \frac{\d k}{ik}.
  \end{align*}
  But $r(k |\omega_n|^{-1/n}t^{-1/n})t \goto 0$ as $t \goto 0$.  To see that the limit can be passed inside the integral, deform $C$ so that it passes along the steepest descent paths of $e^{-i \omega_n k}$, then pass the limit inside using the dominated convergence theorem and deform back to $C$.  From this the result follows for the case of one discontinuity, with compact support.  The general case follows from Theorem~\ref{Thm:Main} below.
\end{proof}

\section{Gibbs phenomenon for dispersive PDEs}
\label{s:generaldispersion}

We now discuss the implications of decomposition~\eref{e:qdecomp}
regarding the behavior of the solution of the IVP in the short-time limit.  We have seen that, apart from a constant offset, 
the dominant behavior of the solution in the regularization region near a discontinuity of the IC is provided 
by the function $I_{\omega,0}(y,t)$.
In this section we therefore examine more closely the properties of such functions.
We start by discussing a simple example.

\subsection{Example: Heat equation.}

Consider the PDE
\[
q_t = q_{xx}\,,
\label{e:heat}
\]
corresponding to $\omega(k) = - ik^2$.
Let $s=y/t^{1/2}$ and $\lambda= kt^{1/2}$. Then
\[
I_{\heat,0}(y,t) = \frac1{2\pi}\int_C \e^{i\lambda s - \lambda^2}\frac{\d\lambda}{i\lambda}
 = \half\,(\erf(s/2)-1)\,,
\]
where with some abuse of notation we write $I_{\heat,0}(y(s),t) = I_{\heat,0}(s)$.
Note that an easy way to compute the above integral is by using the relation
\[
\deriv{ }s I_{\heat,0}(s) = \frac1{2\pi}\int_C\e^{i\lambda s - \lambda^2}\d\lambda = \frac1{2\sqrt\pi}\e^{-s^2/4}\,,
\label{e:dI0ds:heat}
\]
We will see a generalization of \eref{e:dI0ds:heat} later.

\begin{figure}[t!]
\kern-\bigskipamount
\centerline{\includegraphics[width=0.405\textwidth]{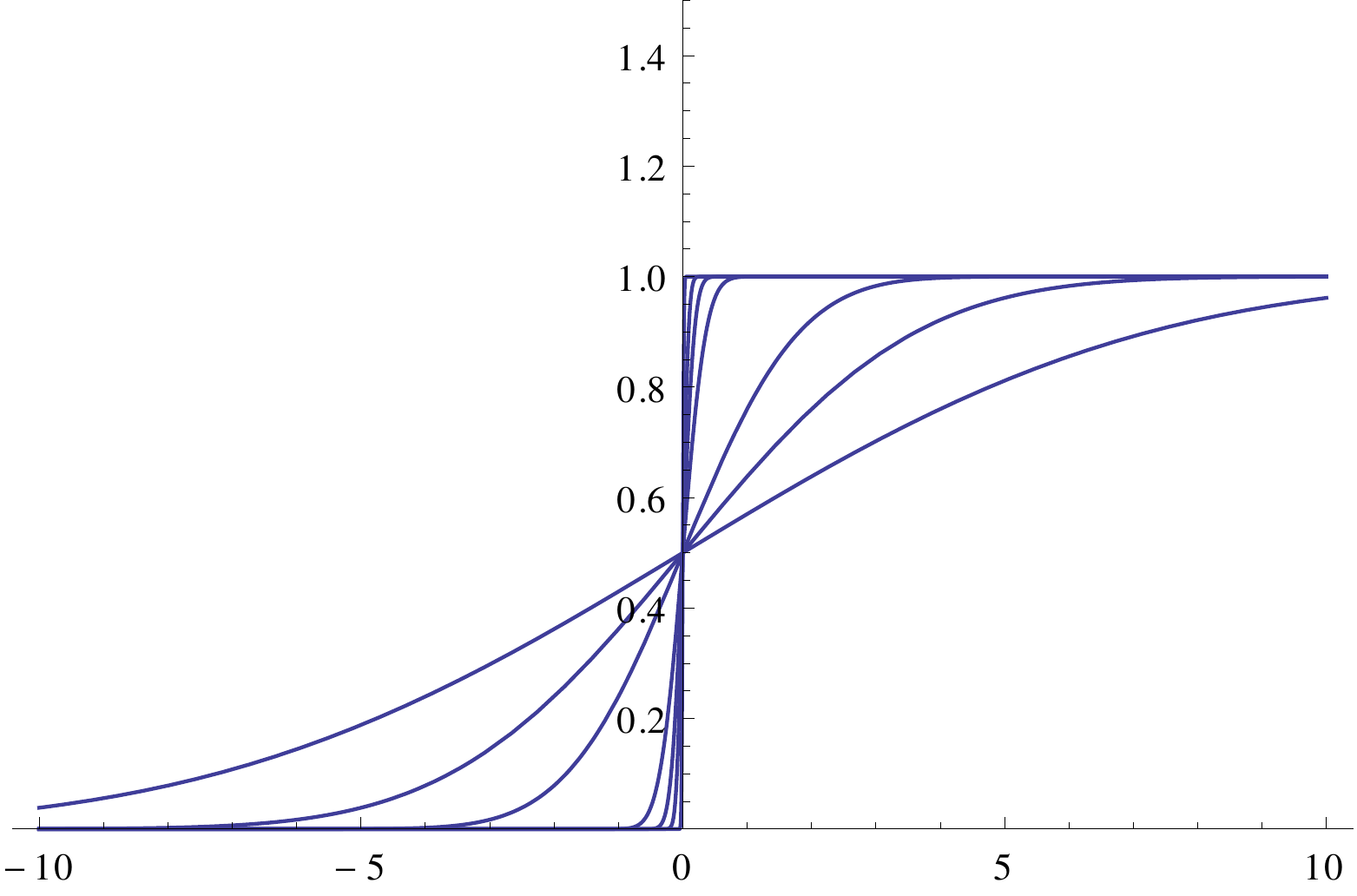}}
\caption{The integral $I_0(x,t)+1$ (vertical axis) as a function of $x$ (horizontal axis) for the heat equation~\eref{e:heat} 
at various values of time: $t=0.01$, 0.05, 0.1, 0.2, 1, 2, and 4.}
\label{f:heat}
\end{figure}

Figure~\ref{f:heat} shows the value of $I_0(x,t)$ as a function of $x$ at different times.
The resulting effect is that of a diffusion-induced smoothing of the initial discontinuity.
This behavior is well-known, and is discussed in most classical PDE books \cite{evans}.
What is perhaps less known, however, is the counterpart of this behavior for dispersive PDEs, which we turn to next.

\subsection{Example: Schr\"odinger equation.}

Consider now the free-particle, one-dimensional linear Schr\"o\-dinger equation, namely, 
\[
iq_t+q_{xx}=0\,,
\label{e:schrodinger}
\]
corresponding to $\omega(k) = k^2$.  In this case,
\[
I_{\schr,0}(s) = \frac1{2\pi i}\int_C \e^{i\lambda s - i\lambda^2}\frac{\d\lambda}\lambda
 = \half\,(\erf(\e^{-i\pi/4}s/2)-1)\,.
\]

\begin{figure}[t!]
\centerline{\includegraphics[width=0.405\textwidth]{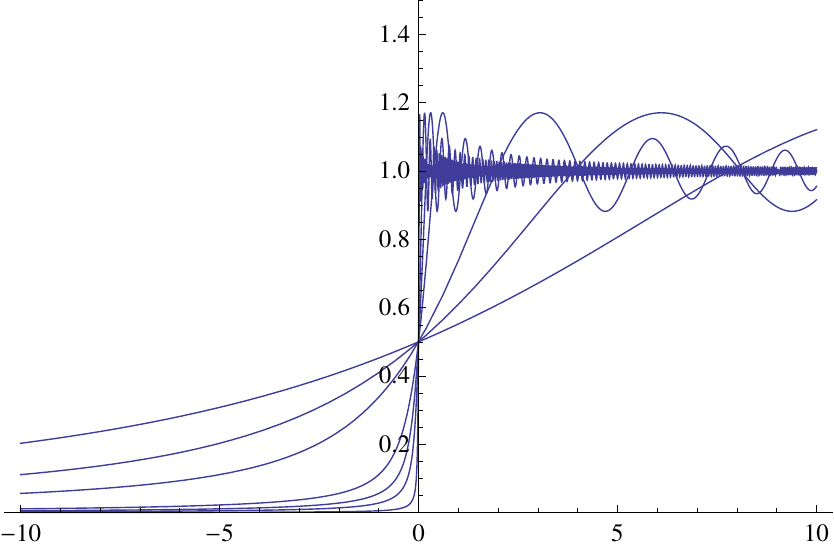}\qquad
\includegraphics[width=0.405\textwidth]{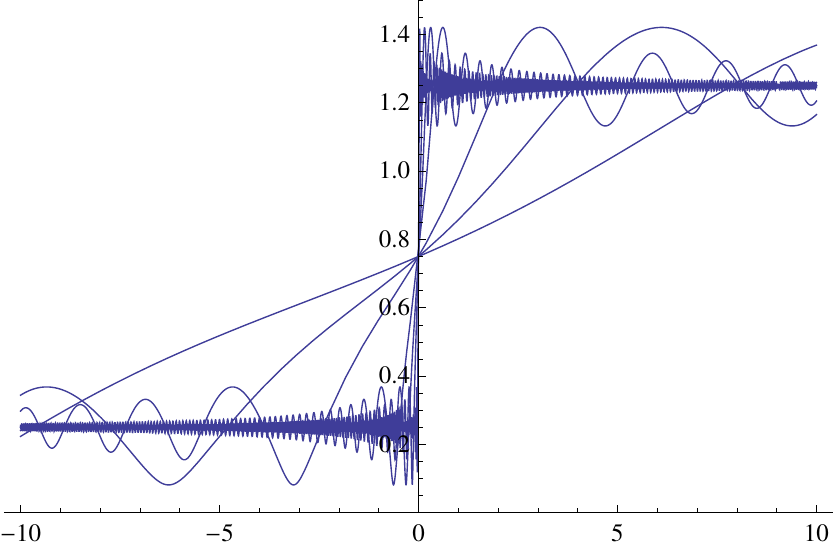}}
\caption{Left: 
Absolute value $|I_0(x,t)+1|$ as a function of $x$ for the Schr\"odinger equation~\eref{e:schrodinger} 
at the same values of $t$ as in Fig.~\protect\ref{f:heat}.
Right: 
Same for $|I_0(x,t)+\txtfrac54|$.
Note in this last case the presence of oscillations to the left of the jump.}
\label{f:schrodinger}
\end{figure}

The corresponding behavior is shown in Fig.~\ref{f:schrodinger}.
For both PDEs, the dominant behavior near the discontinuity is expressed in terms of a \textit{similarity solution},
depending on $x$ and $t$ only through the similarity variable $s=(x-c)/t^{1/2}$, as seen in Theorem~\ref{t:universality}.
The solution behavior however is very different:
While for the heat equation the integral $I_{\omega,0}(x,t)$ captures the smoothing effect of the PDE,
for the Schr\"odinger equation, $I_{\omega,0}(x,t)$ results in oscillations.


\subsection{Example: Stokes equation.}

\begin{figure}[t]
\centerline{\includegraphics[width=0.405\textwidth]{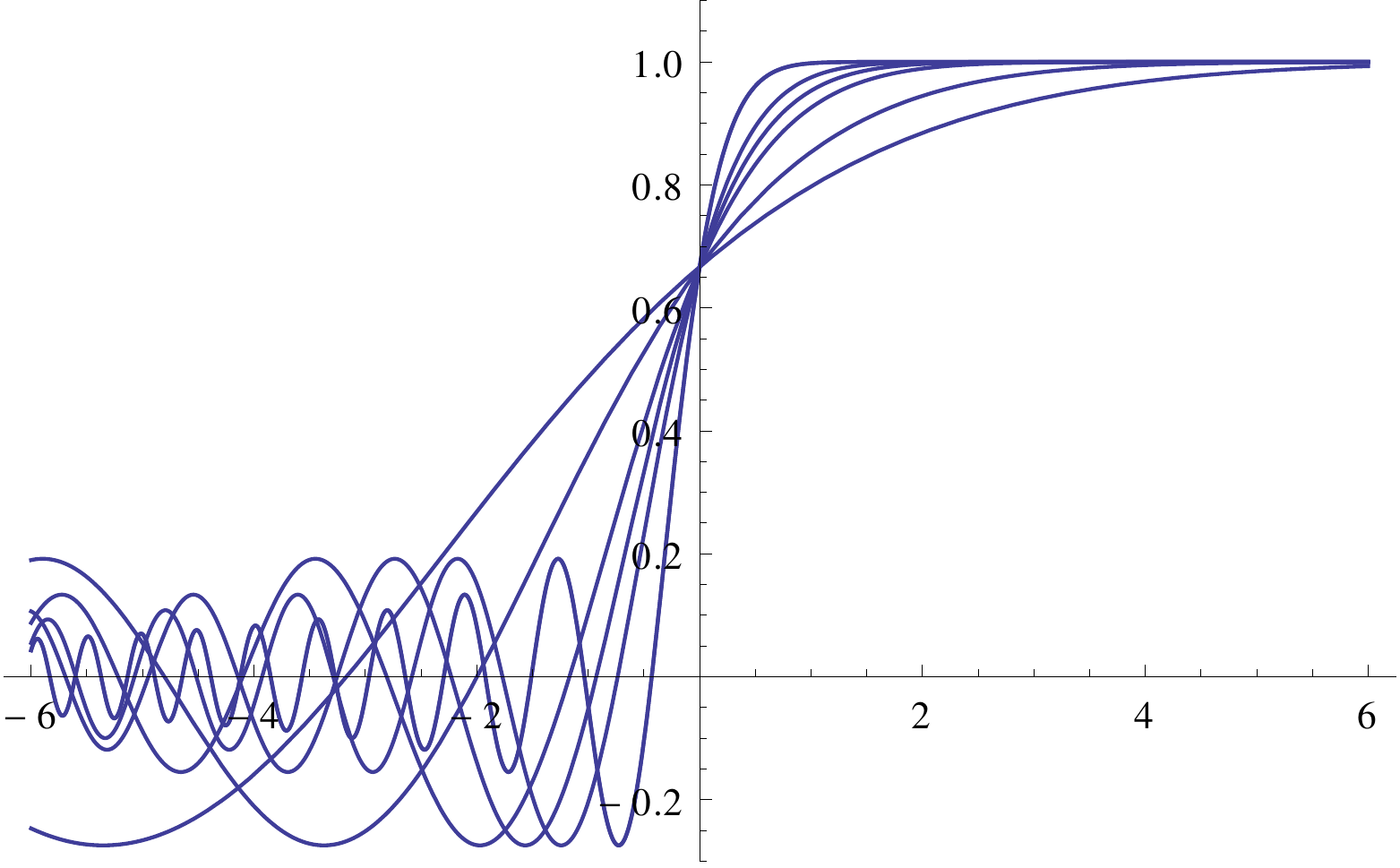}}
\caption{Same as Fig.~\ref{f:heat}, but for the Stokes equation~\eref{e:stokes}.}
\label{f:stokes}
\end{figure}

Consider now the Stokes equations 
\[
q_t + q_{xxx} = 0\,,
\label{e:stokes} 
\]
corresponding to $\omega(k) = -k^3$.  Letting $s= y/t^{1/3}$ and $\lambda = kt^{1/3}$ one has,
using similar methods as before,
\[
I_{\mathrm{stokes},0}(y,t) = \frac1{2\pi i}\int_C \e^{i\lambda s - i\lambda^3}\frac{\d\lambda}\lambda
 = \int_{s/\sqrt[3]3}^\infty \!\Ai(z)\,\d z\,,
\]
where $\Ai(z)$ is  the classical Airy function (e.g., see~\cite{lebedev,NIST2010}), which admits the integral representation
$\Ai(z) = \int\nolimits_\Real \e^{i\lambda z - i\lambda^3}\d\lambda/(2\pi)$ 
\cite{bleisteinhandlesman}.
The corresponding behavior is illustrated in Fig.~\ref{f:stokes}.

Note that, since all the PDEs considered in this work are linear, the behavior arising from a negative jump 
is simply the reflection with respect to the horizontal axis of that for a positive jump.
On the other hand, unlike the heat and Schr\"odinger equation, the Stokes equation does not possess left-right symmetry.
So the values of $I_{\omega,0}(y,t)$ to the left of the discontinuity are not symmetric to those to the right
(as is evident from Fig.~\ref{f:stokes}).
Note also that 
the results for the Stokes equation with the opposite sign of dispersion 
(i.e., $q_t - q_{xxx} = 0$) 
are obtained by simply exchanging $x-c$ with $c-x$ (i.e., $y$ with $-y$) in the above discussion.

\subsection{Gibbs-like oscillations of dispersive PDEs.}\label{s:gibbs}
The solution of the Schr\"o\-dinger equation described above shares the three defining features of the Gibbs phenomenon, namely:
(i) non-uniform convergence of the solution of the PDE to the IC as $t\downarrow0$ in a neighborhood of the discontinuity;
(ii) spatial oscillations with increasing (in fact, unbounded) frequency as $t\downarrow0$ (because they are governed by the similarity variable);
(iii) constant overshoot in a neighborhood of the discontinuity as $t\downarrow0$. (We will elaborate on this last issue later in the section.)
Thus, the limit $t\downarrow0$ for the solution of the PDE is perfectly analogous to the limit $n\to\infty$ in the truncation of the Fourier series.




Recall that, while $q_c$ contributes a constant offset to the solution, 
the value of $q(x,t)$ at $(c,0)$ [as obtained from the reconstruction formula~\eref{e:qdecomp}]
will differ from $q_c$, because, even though $q_\res(0,0)=0$, in general, $\lim_{t \downarrow 0}I_{\omega,0}(0,t)\ne0$.  For monomial dispersion relations, i.e., $\omega_n(k) = \omega_nk^n$, 
it easy to see that $I_{\omega_n,0}(0,t)$ is actually independent of time.
In fact, the value of $I_{\omega_n,0}(0,t)$ can be easily obtained explicitly.
From~\eref{e:I0omegadef} we have
\begin{equation*}
I_{\omega_n,0}(0,t) 
  = \frac1{2\pi} \pvint_{\mathbb R} \e^{\pm i\lambda^n}\frac{\d\lambda}{i\lambda}
    - \frac12\,,
\end{equation*}
since $\Res\nolimits_{\lambda=0}[\e^{\pm i\lambda^n}/(i\lambda)]=1$.
Now note that 
$\pvint\nolimits_{\mathbb R} \e^{\pm i\lambda^n}\d\lambda/(i\lambda)=0$ for $n$ even,
while the same integral equals $\pm\int\nolimits_{\mathbb R}\sin(\lambda^n)\,\d\lambda/\lambda = \pm\pi/n$ for $n$ odd.
Hence we have simply
\begin{equation}
I_{\omega_n,0}(0,t) = \begin{cases}
  -\frac12\,, &n~\mathrm{even}\,,\\
  -\frac12(1\pm 1/n)\,, & n~\mathrm{odd}\,.
\end{cases}
\end{equation}

One can carry out the analogy with the classical Gibbs phenomenon even further 
and compute the ``overshoot'' of these special functions --- namely, the ratio of 
the maximum difference between the value of the special function and the jump, compared to the jump size.
Recall that the overshoot for the Gibbs phenomenon is given by the Wilbraham-Gibbs constant \cite{gibbs,wilbraham} (see also \cite{hewitt}).
\[
\mathfrak g = \frac{1}{\pi} \int_0^\pi \frac{\sin z}{z} \,\d z  - \half \approx 0.089490\ldots
\label{e:gibbsconstant}
\]
For example, the maximum value of the partial sum of the Fourier series for $\chi_{[-1,1]}(y)$ on $[-2,2]$ 
will converge to $1+\mathfrak g$, and its minimum to $-\mathfrak g$.  

To examine the overshoot of the special functions, 
we look at $G_n(y,t) = I_{\omega_n,0}(y,t) + 1$, which converges pointwise to $\chi_{(0,\infty)}(y)$ for all $y \ne 0$ as $t \downarrow 0$.
Specifically, we compute numerically the maximum and minimum of the real part, imaginary part and modulus of $G_n(y,t)$.
Note that, for all $t\ne0$, all such values are independent of $t$.   
Table~\ref{f:table} shows these values as a function of $n$.
Surprisingly, the table shows that 
\textit{these values converge to exactly the same constants as for the Gibbs phenomenon as $n\to\infty$!}

\begin{table}[t!]
\small
\begin{tabular}{c|cc|cc|c}
$n$ & max real & min real & max imag & min imag & max modulus \\
\hline
2 & 1.17025 & -0.170246 & 0.243797 & -0.243797 & 1.17066 \\
3 & 1.27435 & 0 & 0 & 0 & 1.27435 \\
4 & 1.11501 & -0.115008 & 0.121603 & -0.121603 & 1.10603 \\
5 & 1.19824 & -0.0159841 & 0 & 0 & 1.19824 \\
6 & 1.10146 & -0.101461 & 0.0819619 & -0.0819619 & 1.10103 \\
7 & 1.16611 & -0.0308676 & 0 & 0 & 1.16611 \\
8 & 1.0963 & -0.0962954 & 0.0618324 & -0.0618324 & 1.09625 \\
9 & 1.14849 & -0.0413221 & 0 & 0 & 1.14849 \\
10 & 1.09384 & -0.0938431 & 0.0496286 & -0.0496286 & 1.09383 \\
11 & 1.1374 & -0.0487894 & 0 & 0 & 1.1374 \\
60 & 1.08961 & -0.0896059 & 0.0083311 & -0.0083311 & 1.08961 \\
120 & 1.08952 & -0.0895187 & 0.00416638 & -0.00416638 & 1.08952 \\
180 & 1.0895 & -0.0895026 & 0.00277769 & -0.00277769 & 1.0895 \\
240 & 1.0895 & -0.089497 & 0.0020833 & -0.0020833 & 1.0895 \\
300 & 1.08949 & -0.0894945 & 0.00166665 & -0.00166665 & 1.08949 
\end{tabular} 
\caption{Numerically computed values for the maximum and minimum of the real part, imaginary part and modulus of $G_n(y,t) = 1 +I_{n,0}(y,t)$ 
as a function of~$n$.  
The overshoot converges to the Wilbraham-Gibbs constant $\mathfrak g$ [cf.~\eref{e:gibbsconstant}].}
\label{f:table}
\end{table}

Indeed, a simple calculation shows why this is true.  
Integration by parts or a simple change of variable can be used to show that, as $n\to\infty$, 
\begin{align*}
I_{\omega_n,0}(y,1) = \frac{1}{2\pi i} \int_{C'} e^{iky - i k^n} \frac{\d k}{k} + \mathcal O(1/n)\,,
\end{align*}
where $C'= C\cap \{k\in\Complex:|\real k| \leq 1\}$, and where without loss of generality 
the semi-circle component of $C$ was taken to have radius less than one. 
Then, by the dominated convergence theorem
we have
\vspace*{-0.4ex}
\begin{align*}
\lim_{n \goto \infty} \frac{1}{2\pi i} \int_{C'} e^{iky - i k^n} \frac{\d k}{k} = \frac{1}{2\pi i} \int_{C'} e^{iky} \frac{\d k}{k}\,,
\end{align*}
where convergence is uniform in $y$. Moreover, the integral on the RHS is easily shown to be
\vspace*{-0.4ex}
\begin{align*}
\frac{1}{2\pi i} \int_{C'} e^{iky} \frac{\d k}{k} = \frac{1}{2\pi} \int_{-1}^1 \frac{\sin ky}{k} \d k - \half\,,
\end{align*}
where the contour in the RHS was deformed back to the real axis since there is a removable singularity at $k=0$.
After a simple rescaling we then have
\vspace*{-0.4ex}
\[
\lim_{n\to\infty} I_{\omega_n,0}(y,1) = \frac1\pi\int_0^\pi \frac{\sin(\pi yz)}z\,\d z - \half\,,
\]
uniformly in $y$.
This integral is maximized and minimized at $y= \pm1$, respectively, yielding
\vspace*{-0.4ex}
\begin{gather*}
\lim_{n \goto \infty} \sup_{y \in \mathbb R} \real G_n(y,1) = 1+\mathfrak g\,,\qquad 
\lim_{n \goto \infty} \sup_{y \in \mathbb R} \imag G_n(y,1) = 0\,,\\
\lim_{n \goto \infty} \inf_{y \in \mathbb R} \real G_n(y,1) = -\mathfrak g\,,\qquad
\lim_{n \goto \infty} \inf_{y \in \mathbb R} \imag G_n(y,1) = 0\,.
\end{gather*}
Note that for a fixed value of $n$ such maxima and minima can occur on either side of the jump
(e.g., cf.\ Figs.~\ref{f:schrodinger} and~\ref{f:stokes}).

\begin{proof}[Proof of Theorem~\ref{t:gibbs}]
  The solution $q(x,t)$ is given by
  \begin{align*}
    q(x,t) = I_{\omega,0}(x+c,t) - I_{\omega,0}(x-c,t), ~~~ \omega(k) = k^n.
  \end{align*}
  Near $x = -c$ we have
  \begin{align*}
    q(y-c,t) = G_n(y,t) - (I_{\omega,0}(y-2c,t) + 1),\quad y \in (-\delta,\delta), ~~ 0 < \delta < 2c.
  \end{align*}
  It follows from Theorem~\ref{Thm:Kernel} that
  \begin{align*}
    |I_{\omega,0}(y-2c,t) + 1| \leq C_\delta t^{1/(2n-2)}, \quad C_\delta > 0,
  \end{align*}
  uniformly for all $y \in (-\infty,\delta)$.  So,
  \begin{align*}
    \lim_{t\downarrow 0}\left( \sup_{|y| \leq \delta} \real G_n(y,t) - C_\delta t^{1/(2n-2)} \right)  \leq  \lim_{t\downarrow 0} \sup_{|y| \leq \delta} \real q(y-c,t)  \leq \lim_{t\downarrow 0} \left(\sup_{|y| \leq \delta} \real G_n(y,t) + C_\delta t^{1/(2n-2)}\right),
  \end{align*}
  and $\lim_{t\downarrow 0} \sup_{|y| \leq \delta} \real q(y-c,t)  = \sup_{y \in \mathbb R} \real G_n(y,1)$.  From this the first claim in the theorem follows for $\delta < 2 c$.   To allow $\delta$ to be larger, just break the analysis into an interval contained in $(-\infty,0]$ and another interval contained in $[0,\infty)$.  The other claims follow from similar calculations.
\end{proof}

\section{Short-time asymptotics: ICs with discontinuous derivatives}
\label{s:discontinuousderivatives}

We now treat the case where one of the derivatives of $q_o$ is discontinuous.
We begin by assuming a discontinuity in the first derivative, then we treat the general case.
We will further generalize the results in Section~\ref{s:multiplediscontinuities}.

\begin{assume}
Let
\begin{itemize}
\item$q_o \in H^1(\mathbb R)$,
\item$[q'_o(c)] = q_o'(c^+)-q_o'(c^-) \neq 0$,
\item$q''_o$ exists on $(-\infty,c)\cup (c,\infty)$,
\item$q''_o \in L^q(-\infty,c) \cap L^q(c,\infty)$ for some $1 < q < \infty$, and
\item $q_o$ is compactly supported.
\end{itemize}
\end{assume}

Assuming compact support avoids possible complications arising from the non-existence of some principal value integrals. 
(This assumption will be removed in Section~\ref{s:multiplediscontinuities}.)
We will show that the asymptotic behavior in the regularization region is given by integrals of the special 
functions considered in the previous section.

Note first that, if $F(k)$ is analytic in a neighborhood of the origin,
\eref{e:qdecomp0} can be written as
\begin{align*}
q(x,t) = [q_0(c)] I_{\omega,0}(x-c,t)+ \frac{1}{2 \pi} \int_{C} e^{i\theta(x,t,k)} F(k) \frac{\d k}{ik},
\end{align*}
with $I_{\omega,0}(y,t)$ and $F(k)$ given by~\eref{e:I0def} and~\eref{e:Fdef}, respectively,
and with $C$ as in Figure~\ref{f:contour}.  
Analyticity of $F$ is always guaranteed if $q_o$ has compact support.  
In the case that $q_o$ is continuous but $q_o'$ is discontinuous, 
we perform one more integration by parts and write
\begin{gather}
q(x,t) = [q_o'(c)] I_{\omega,1}(x-c,t)+ \frac{1}{2 \pi} \int_{C} e^{i\theta(x,t,k)} F_1(k) \frac{\d k}{(ik)^2}\,,
\label{e:C-form}
\\
\noalign{\noindent where}
F_1(k) = \left( \int_{-\infty}^c + \int_{c}^\infty \right) e^{-iks} q_o''(s) \d s,
\notag
\end{gather}
and where we have introduced the generalization of $I_{\omega,0}(y,t)$ as
\begin{gather}
I_{\omega,m}(y,t) = \frac{1}{2\pi} \int_{C} \frac{e^{i ky - i \omega(k) t}}{(ik)^{m+1}} \d k\,.
\label{e:Iomegandef}
\end{gather}

As before, we now expand~\eref{e:C-form} both near and away from the singularity $c$.
In a neighborhood of $(c,0)$, we leave $I_{\omega,1}(y,t)$ alone, and we expand $F_1(k)$.  
As $k\to0$, 
\begin{align*}
e^{ikc} e^{i\theta(x-c,t,k)} = e^{ikc} ( 1 + ik(x-c) + \bigo(k^2))\,.
\end{align*}
We then have
\begin{align*}
q(x,t) = [q_o'(c)] I_{\omega,1}(x-c,t)+ \frac{1}{2 \pi} \int_{C} e^{ikc} \left(\frac{ 1 + ik(x-c)}{(ik)^2} \right) F_1(k) \,\d k
+ q_{\res,1}(x-c,t),
\end{align*}
where
\begin{align*}
q_{\res,1}(x-c,t) = \frac{1}{2 \pi} \int_{\mathbb R} e^{ikc} \left(\frac{e^{i\theta(x-c,t,k)}- 1 - ik(x-c)}{(ik)^2} \right) F_1(k) \,\d k\,.
\end{align*}
We expect $q_{\res,1}(y,t)$ to give a lower order contribution as $(x,t) \rightarrow (c,0)$.
We thus examine this expression in the regularization region $|x-c| \leq C t^n$. 
Lemma~\ref{Lemma:ResEst-HO} indicates that $q_{\res,1}(x,t) = \bigo(t^{1/n+1/(np)})$ 
because $F \in L^p(\mathbb R)$ (where $1/p + 1/q = 1$). 
Therefore $q_{\res,1}(y,t)$ can indeed be seen as the error term.

We now examine~\eqref{e:C-form} for $|x-c| \ge \delta > 0$ and $|s-x| \le \delta/2$.
We have
\begin{align*}
q(x,t) - q_o(s)  &= [q_o'(c)] (I_{\omega,1}(x-c,t)-I_{\omega,1}(s-c,0))
+ \frac{1}{2 \pi} \int_{C}e^{iks} (e^{i\theta(x-s,t,k)}-1) F_1(k) \frac{\d k}{(ik)^2}
\\
&= 
[q_o'(c)] (I_{\omega,1}(x-c,t)-I_{\omega,1}(s-c,0))
+ \frac{(x-s)}{2 \pi} \int_C e^{iks} F_1(k) \frac{\d k}{ik} + q_{\res,1}(x-s,t).
\end{align*}
Applying Theorem~\ref{Thm:Kernel} and Lemma~\ref{Lemma:ResEst-HO}, in the regularization region $|x-s|^n \le Ct$ we have
\begin{align*}
q(x,t) &= q_o(s) + [q_o'(c)] ( (s-c)\chi_{(-\infty,c)}(s) - (x-c)\chi_{(-\infty,c)}(x))\\
&+ \frac{(x-s)}{2 \pi} \int_C e^{iks} F_1(k) \frac{\d k}{ik} + \bigo\left( t^{3/(2(n-1))} + t^{1/n+ 1/(np)} \right).
\end{align*}
This expression is simplified using $\chi_{(-\infty,c)}(s) = \chi_{(-\infty,c)}(x)$ and the relation
\begin{align*}
\frac{(x-s)}{2 \pi} \int_C e^{iks} F_1(k) \frac{\d k}{ik} = -\half [q_o'(c)](x-s) + \frac{(x-s)}{2 \pi} \pvint_C e^{iks} F_1(k) \frac{\d k}{ik},
\end{align*}
to obtain
\begin{align*}
q(x,t) &= q_o(s) + [q_o'(c)](s-x)(-1/2 + \chi_{(-\infty,c)}(s))
\\
& + \frac{(x-s)}{2 \pi} \pvint_C e^{iks} F_1(k) \frac{\d k}{ik} 
  + \bigo\left( t^{3/(2(n-1))} + t^{1/n+ 1/(np)} \right).
\end{align*}

Next we generalize the above result to a discontinuity in a derivative of arbitrary order:
\begin{assume} 
\label{Assume:m}
Let
\begin{itemize}
\item$q_o \in H^m(\mathbb R)$,
\item$[q^{(m)}_o(c)] \neq 0$,
\item$q^{(m+1)}_o$ exists on $(-\infty,c) \cup (c,\infty)$, separately,
\item$q^{(m+1)}_o \in L^q(-\infty,c) \cap L^q(c,\infty)$ for some $1 < q < \infty$, and
\item $q_o$ is compactly supported.
\end{itemize}
\end{assume}
Let $a_\ell(y,t)$ be the Taylor coefficients of $e^{i\theta(y,t,k)}$ at $k = 0$.  
Then for $s \in \mathbb R$ (possibly equal to $c$) we find the expansion
\begin{multline}\label{e:ho-expand}
q(x,t) = [q_o^{(m)}(c)] I_{\omega,m}(x-c,t) + \frac{1}{2 \pi} \int_{C} e^{iks} \left( \sum_{\ell=0}^m a_\ell(x-s,t) k^\ell \right) F_m(k) \frac{\d k}{(ik)^{m+1}} + q_{\res,m}(x-s,t),
\end{multline}
where
\begin{align*}
q_{\res,m}(x-s,t) &= \frac{1}{2 \pi} \int_{\mathbb R} e^{iks} \left(e^{i\theta(x-s,t,k)}-\sum_{\ell=0}^m a_\ell(x-s,t) k^\ell \right) F_{m}(k) \frac{\d k}{(ik)^{m+1}},\\
F_m(k) &= \left( \int_{-\infty}^c + \int_{c}^\infty \right) e^{-ikx} q_o^{(m+1)}(x) \d x.\notag
\end{align*}
Invoking Lemma~\ref{Lemma:ResEst-HO}, this expression provides the asymptotic expansion in the regularization region 
$|x-s|^n \leq Ct$. 
Indeed, $q_{\res,m}(x,t) = \bigo( t^{m/n+ 1/(pn)} )$ for $1/p + 1/q = 1$.  
This expansion can be understood more thoroughly as follows.   
Formally, for $s \in \mathbb R$
\begin{multline}
\label{e:taylorterms}
(-i\partial_x)^j q_o(s) = [q_o^{(m)}(c)] \Res_{k=0} \left( \frac{e^{ik(s-c)}}{i (ik)^{m-j+1}} \right)\chi_{(-\infty,0)}(s-c) + \frac{1}{2\pi} \int_{C} e^{iks} F_{m}(k) \frac{\d k}{(ik)^{m-j+1}}.
\end{multline}
We next show that
\begin{align}\label{series-approx}
\sum_{j=0}^M \frac{(-it)^j}{j!} \omega(k)^j = \sum_{\ell = 0}^{nM} a_\ell(0,t) k^\ell + \bigo(t^{M+1}k^{nM}),
\end{align}
as $|k| \goto \infty$ and $t \downarrow 0$. 
To see this, it follows from Lemma~\ref{Lemma:ResEst-HO} that $a_\ell(0,t) =\mathcal O(t^{\ell/n})$ and then
\begin{align*}
e^{-i\omega(k)t}-\sum_{\ell=0}^{nM} a_\ell(0,t) k^\ell  = \bigo(t^{M+1})
\end{align*}
as $t \downarrow 0$, because only integer powers of $t$ appear.  Then 
\begin{align*}
e^{-i\omega(k)t}-\sum_{j=0}^M \frac{(-it)^j}{j!} \omega(-i\partial_x)^j  = \bigo(t^{M+1}),
\end{align*}
implying
\begin{align*}
\sum_{j=0}^M \frac{(-it)^j}{j!} \omega(-i\partial_x)^j = \sum_{\ell = 0}^{nM} a_\ell(0,t) k^\ell + \bigo(t^{M+1}).
\end{align*}
Then \eqref{series-approx} follows by noting that both sides have no powers of $k$ larger than $k^{nM}$.  
In turn, \eqref{series-approx} implies
\begin{align}\label{e:ic-expand}
\sum_{j=0}^M \frac{(-it)^j}{j!} \omega(-i\partial_x)^j q_o(s) &= [q_o^{(m)}(c)] \Res_{k=0} \left(\frac{e^{ik(s-c)}}{i(ik)^{m+1}} \sum_{\ell=0}^{nM} a_\ell(0,t)k^\ell \right) \chi_{(-\infty,0)}(s-c)\\
&+ \frac{1}{2\pi} \int_{C} e^{iks} \left(  \sum_{\ell=0}^{nM} a_\ell(0,t)k^\ell\right) F_{m}(k) \frac{\d k}{(ik)^{m+1}} + \bigo(t^{M+1}).
\end{align}
If $s \neq c$ then this expression is well-defined and continuous for $nM \leq m$.  
If $s =c$, there are issues concerning the definition of the value of $q^{(nM)}_o(c)$ on the left-hand side of the equation and we must restrict to $nM < m$.

\paragraph{Near the singularity.}

Let $M = \lfloor (m-1)/n \rfloor$.  For $|x-c|^n \leq C t$ we combine \eqref{e:ic-expand} and \eqref{e:ho-expand}  to find
\begin{align}\label{e:in}
\begin{split}
q(x,t) &= \sum_{j=0}^M \frac{(-it)^j}{j!}\omega(-i\partial_x)^j q_o(c)  +  [q_o^{(m)}(c)] I_{\omega,m}(x-c,t)\\
&+ \frac{1}{2 \pi} \int_{C} e^{ikc} \left( \sum_{\ell=0}^m ( a_\ell(x-c,t) -a_\ell(0,t)) k^\ell \right) F_m(k) \frac{\d k}{(ik)^{m+1}} + \bigo \left(t^{\frac{m}{n}+ \frac{1}{np}} \right).
\end{split}
\end{align}
Here, the residue term in \eref{e:taylorterms} vanishes at $s=c$ because $Mn < m$ and no $k^{-1}$ term is present. It also follows (see Lemma~\ref{Lemma:ResEst-HO}) that $a_\ell(x-c,t) = \bigo(t^{\ell/n})$ so that this is indeed a consistent expansion.

 \paragraph{Away from the singularity.}

Let $M = \lfloor m/n \rfloor$.  We examine the expansion for  near $x = s$ for $|s-c| \geq \delta  >0$.  We use the short-time asymptotics for $I_{\omega,m}$ (see Theorem~\ref{Thm:Kernel}) to find for $|x-s|^n \leq C |t|$
\begin{align}\label{e:out}
\begin{split}
q(x,t) &= \sum_{j=0}^M \frac{(-it)^j}{j!}\omega(-i\partial_x)^j q_o(s)\\
&+ \frac{1}{2 \pi} \int_{C} e^{iks} \left( \sum_{\ell=0}^m ( a_\ell(x-s,t) -a_\ell(0,t)) k^\ell \right) F_m(k) \frac{\d k}{(ik)^{m+1}}\\
& - i[q_o^{(m)}(c)] \Res_{k=0} \left( \frac{e^{ik(x-c)-i \omega(k) t}}{(ik)^{m+1}} - \frac{e^{ik(s-c)}}{(ik)^{m+1}} \sum_{j=0}^M \frac{(- i \omega(k)t)^j}{j!} \right) \chi_{(-\infty,0)}(s-c) \\
&+ \bigo \left(t^{\frac{m}{n} } \left( t^{\frac{1}{np}} + t^{\frac{n+2m}{2n(n-1)}} \right) \right).
\end{split}
\end{align}
If we set $x = s$ then the residue term is $\bigo(t^{M+1})$ ($m/n+ 1/n \leq M+1$) and the short-time Taylor expansion 
\begin{align}\label{e:Taylor}
q(x,t) &= \sum_{j=0}^M \frac{(-it)^j}{j!} \omega(-i\partial_x)^j q_0(x) + \bigo \left(t^{\frac{m}{n} } \left( t^{\frac{1}{np}} + t^{\frac{n+2m}{2n(n-1)}} \right) \right).
\end{align}
follows.  Here the error term is uniform in $x$ as $x$ varies in the region $|x-c| \geq \delta$.  Thus, in particular, if $q_o$ vanishes identically in a neighborhood of $s$ then for $|x-s|^n \leq Ct$
\begin{align}
q(x,t) =  \bigo \left(t^{\frac{m}{n} } \left( t^{\frac{1}{np}} + t^{\frac{n+2m}{2n(n-1)}} \right) \right).\label{e:ZeroTaylor}
\end{align}

\paragraph{A unified formula.}

We now introduce some convenient and unifying notation that will be useful to combine the above results.  
Define
\begin{align*}
R_{M,m,c}(q_o;x,s) &= -i[q_o^{(m)}(c)] \Res_{k=0} \left( \frac{e^{ik(x-c)-i \omega(k) t}}{(ik)^{m+1}} - \frac{e^{ik(s-c)}}{(ik)^{m+1}} \sum_{j=0}^M \frac{(- i \omega(k)t)^j}{j!} \right) \chi_{(-\infty,0)}(s-c),\\
A_{m}(q_o;x,s) &= \frac{1}{2 \pi} \int_{C} e^{iks} \left( \sum_{\ell=0}^m ( a_\ell(x-s,t) -a_\ell(0,t)) k^\ell \right) F_m(k) \frac{\d k}{(ik)^{m+1}}\,.
\end{align*}
Note $A_m(q;x,s)$ can only be applied to functions whose Fourier transform is analytic in a neighborhood of the origin. 
Therefore we have for $M = 0, \ldots, \lfloor \frac{m-1}{n} \rfloor$ and $s \in \mathbb R$,
\begin{align*}
q(x,t) = \sum_{j=0}^M \frac{(-it)^j}{j!} \omega(-i\partial_x)q_o(s) + A_{m}(q_o;x,s)
+ \begin{cases} R_{M,m,c}(q_o;x,s),& s \neq c,\\
\\
[q_o^{(m)}(c)] I_{\omega,m}(x-c,t),& s = c,
\end{cases}
+ \bigo \left(t^{\frac{m}{n} } \left( t^{\frac{1}{np}} + t^{\frac{n+2m}{2n(n-1)}} \right) \right).
\end{align*}
While the formula for $s \neq c$ is also valid for $M = \lfloor m/n \rfloor$, this is a convenient form.  
Furthermore, when no singularity is present locally, \eref{e:Taylor} provides a cleaner formula in terms of quantities that are easier to compute.  
We note that $A_m$ and $R_{M,m,c}$ ($s \neq c$) contain terms that are analytic in $x$ and $t$ while $I_{\omega,m}$ encodes the 
dominant behavior near the singularity; \emph{i.e.}, it has a discontinuous derivative at some order.

\section{Short-time asymptotics: ICs with multiple singular points and non-compact support}
\label{s:multiplediscontinuities}

We now discuss the case of ICs with multiple points of discontinuity.  
The results in this section are the most general ones of this work regarding the short-time behavior of the solution of dispersive PDEs.

\begin{assume}
\label{Assume:General}
For $c_0 = -\infty < c_1 < \cdots < c_N < c_{N + 1} = +\infty$, let
\begin{itemize}
\item$q_o \in H^m(\mathbb R) \cap L^1((1+|x|)^\ell\d x)$, with $\ell \geq \mathfrak C_n$,
\item$[q^{(m)}_o(c_i)] \neq 0$\, for\, $i = 1,\ldots,N$,
\item$q^{(m+1)}_o(x)$ exists on $(c_{i-1},c_i)$\, for\, $i = 1,\ldots,N+1$,
\item$q^{(m+1)}_o \in  L^2(c_{i-1},c_i)$\, for\, $i = 1,\ldots,N+1$.
\end{itemize}
\end{assume}
Note that we have removed the assumption of compact support.  
The key to do so is to use a Van der Corput \textit{neutralizer} (or ``bump'' function)  (e.g., see \cite{bleisteinhandlesman}),
namely a function that interpolates infinitely smoothly between 0 and 1.
More precisely, for our purposes a neutralizer is a function $\eta_\delta(y)$ with the following properties:
\begin{itemize}
\item[(i)] it possesses continuous derivatives of all orders;
\item[(ii)] $\eta_\delta(y)=1$ for $y<\delta/2$ and $\eta_\delta(y)=0$ for $y>\delta$;
\item[(iii)] the derivatives of $\eta_\delta(y)$ of all orders vanish at $y=\delta/2$ and $y=\delta$.
\end{itemize}
A suitable definition is given by
\[
\eta_\delta(y) = n(\delta-x)/[n(y-\delta/2)+n(\delta-x)]\,
\nonumber
\]
where
\[
n(y) = \begin{cases}
  1, &y<0\,,\\
  \e^{-1/y}, &y>0\,,
\end{cases}
\nonumber
\]
but the actual form of the neutralizer is irrelevant for what follows.
Then, to study the behavior near each discontinuity $(x,t) = (c_j,0)$, for $j=1,\dots,N$, one can decompose the IC as
\begin{gather}
q_o(x) = \sum_{j=1}^m q_{o,j}(x) + q_{o,\reg}(x)\,,
\\
\noalign{\noindent
where}
q_{o,j}(x) = q_o(x)\,\eta_\delta(|x-c_j|)\,,
\\
\noalign{\noindent and}
q_{o,\reg}(x) = q_o(x)\,\bigg(1-\sum_{j=1}^m\eta_\delta(|x-c_j|)\bigg)\,,
\end{gather}
with $\delta<\mathrm{min}_{j=1,\dots,m-1}(c_{j+1}-c_j)/2$.
Correspondingly, the solution of the PDE is decomposed as 
\[
q(x,t) = \sum_{j=1}^m q_j(x,t) + q_\reg(x,t)\,.
\]
Note that each $q^{(m)}_{o,j}(x)$ is discontinuous but compactly supported,
while $q_{o,\reg}^{(m)}(x)$ is non-compactly supported but continuous.  
Moreover, $q_{o,j}(c_{j'}) = 0$ for all $j'\ne j$, and $q_{o,\reg}(c_j)=0$ for $j=1,\dots,m$.  
Importantly, it follows that  $q_{o,\reg} \in H^{m+1}(\mathbb R)$. 
Noting that $[q_{o,\reg}^{(m)}(c)]= 0$, with $nM \le m < n(M+1)$,
by~\eref{e:Taylor} we have 
\begin{align*}
q_\reg(x,t) = \sum_{j=0}^M \frac{(-it)^j}{j!} \omega(-i \partial_x)^j q_{o,\reg}(x) + \bigo(t^{m/n+1/(2n)})\,.
\end{align*}
In the regularization region $|x-c_j|^n \leq C t$, all derivatives of $q_{o,\reg}$ vanish identically so that $q_{\reg}(x,t) =  \bigo(t^{m/n+1/(2n)}) = q_{j'}(x,t)$ for $j' \neq j$, see \eqref{e:ZeroTaylor}.

We state our main asymptotic result as a theorem. 
\begin{theorem}\label{Thm:Main}
Suppose Assumption~\ref{Assume:General} holds.
\begin{itemize}
\item If $|x-c_j|^n \leq C|t|$ then for $M = \lfloor \frac{m-1}{n} \rfloor$
\begin{align}\label{e:full}
\begin{split}
q(x,t) &= \sum_{j=0}^M \frac{(-it)^j}{j!} \omega(-i\partial_x)^j q_0(c_j) + [q_o^{(m)}(c_j)] I_{\omega,m}(x-c_j,t) \\
&+ A_m(q_{o,j};x,c_j) +\bigo\left( t^{\frac{m}{n}}\left(t^{\frac{1}{2n}}+ t^{\frac{n+2m}{2n(n-1)}}\right)\right).
\end{split}
\end{align}

\item If $|c_j-x| \geq\delta > 0$ for all $j$ then for $M = \lfloor \frac{m}{n} \rfloor$
\begin{align*}
q(x,t) &= \sum_{j=0}^M \frac{(-it)^j}{j!} \omega(-i\partial_x)^j q_0(x) +\bigo\left( t^{\frac{m}{n}}\left(t^{\frac{1}{2n}}+ t^{\frac{n+2m}{2n(n-1)}}\right)\right).
\end{align*}

\end{itemize}
\begin{proof}
We use linearity. As discussed, we apply \eref{e:Taylor} and \eqref{e:ZeroTaylor} so that $q_{\reg}(x,t) = \bigo(t^{m/n+1/(2n)})$.  The first claim follows from \eref{e:in} and \eref{e:out}. The final claim follows from \eref{e:Taylor}.
\end{proof}
\end{theorem}

From \eqref{e:full} we conclude that near a singularity $q(x,t)$ can be written as $I_{\omega,m}$  plus lower-order and analytic terms.  We not only have an asymptotic expansion but an expansion that separates regularity properly. Furthermore, the expansion about $c_j$ depends only on local properties of $q_o$ through $q_{o,j}$.


\section{Further analysis and computation of the special functions}
\label{s:specialfunctions}

It should be abundantly clear from Sections~\ref{s:onediscontinuity}--\ref{s:multiplediscontinuities} 
that the integrals $I_{\omega,m}(y,t)$ [defined in~\eref{e:Iomegandef}]
play a crucial role in the analysis.  
The detailed properties of these integrals are discussed in Appendix~\ref{a:regularity}.
Here we mention some further properties of these objects and we outline an efficient computational approach
for their numerical evaluation.

\paragraph{Monomial dispersion relations.}

Recall the definition~\eref{e:I0def} of $I_{\omega,0}(y,t)$.
and let $\omega(k) = \omega_n\,k^n$.  Performing the change of variable 
\[
s=y/(|\omega_n|t)^{1/n}\,,\qquad
\lambda= (|\omega_n|t)^{1/n}k\,,
\label{e:k2lambda} 
\]
with some abuse of notation we have that $I_{\omega,0}(y,t)= I_{n,0}^\sigma(y,t)$ is given by
\begin{gather}
I_{n,0}^\sigma(y,t) = E_{n,1}^\sigma(s)\,,
\label{e:I0omegadef}
\\
\noalign{\noindent 
with $\sigma = e^{i \arg(\omega_n)}$, and 
where we have defined}
E_{n,m}^\sigma(s) = \frac1{2\pi}\int_C \e^{i\lambda s -\sigma i\lambda^n}\frac{\d\lambda}{(i\lambda)^m}\,.
\end{gather}
Like their simpler counterparts $I_{\omega,0}(y,t)$, the integrals $I_{\omega,n}(y,t)$ take on a particularly simple form 
in the case of a monomial dispersion relation.  Taking again $\omega_n\in\Real$, we have
\[
I_{n,m}(y,t) = (|\omega_n|t)^{m/n} E_{n,m}^\sigma(s)\,
\]
Now,
\[
E_{n,m}^\mp(s) = \frac1{2\pi} \int_C \e^{i\lambda s \mp i\lambda^n}\frac{d\lambda}{(i\lambda)^{m+1}}\,.
\]
We then have
\[
\deriv{ }s E_{n,m}^\sigma(s) = E_{n,m-1}^\sigma(s)\,.
\label{e:dEds}
\]
So in principle one could obtain $E_{n,m}^\sigma(s)$ by integrating the right-hand side of~\eref{e:dEds}
and by fixing the integration constant appropriately.
In practice, however, it is more convenient to evaluate the integral for $E_{n,m}^\sigma(s)$ directly, 
using the methods discussed below.

\paragraph{General dispersion relations.}

Following arguments from Lemma~\ref{Lemma:Deform}, for $t > 0$, $I_{\omega,m}(y,t)$ 
may be deformed to a contour that is asymptotically on the path of steepest descent for $e^{-i\omega(k)}$.  
Let $C$ be this contour.  
From this deformation, differentiability follows and 
\[
\partial_y^j I_{\omega,m}(y,t) = I_{\omega,m-j}(y,t)\,.
\]
Yet more structure is present. 
A straightforward calculation using integration by parts shows
\begin{align*}
-i t \omega'(-i \partial_y) I_{\omega,m}(y,t) = \frac{1}{\pi} \int_{C} -i t \omega'(k)\frac{e^{i ky - i \omega(k) t}}{(ik)^{m+1}} \d k = y I_{\omega,m}(y,t).
\end{align*}
We thus have obtained the $(n-1)^\mathrm{th}$-order differential equation
\begin{align}\label{e:diff-eq}
\omega'(-i \partial_y) I_{\omega,m}(y,t) =\frac{iy}{t} I_{\omega,m}(y,t),
\end{align}
satisfied by $I_{\omega,m}(y,t)$.


\paragraph{Dissipative PDEs.}

The results in Section~\ref{s:generaldispersion} 
are easily modified when $\omega_n$ is not real,
i.e., when one is dealing with a dissipative PDE.
Recall that, for well-posedness, this can only happen when $n$ is even, in which case $\omega_n = -i|\omega_n|$.

\subsection{Numerical computation of the special functions.}

Next, we discuss the numerical evaluation of $I_{\omega,m}(y,t)$ for all $y$ and $t$.  
First, introduce $\omega_t(k) = \omega(k t^{-1/n})t = \omega_n k^n + \bigo(t^{1/n}k^{n-1})$.  Then
\begin{align*}
I_{\omega,m}(y,t) = t^{(m-1)/n} I_{\omega_t,m}(y t^{-1/n},1).
\end{align*}
It is important that $\omega_t(k) \approx \omega_n k^n$ for $t$ small.  We consider the computation of $I_{\omega_t,m}(s,1)$ accurately for all $s \in \mathbb R$.  The numerical method for accomplishing this follows the proof of Theorem~\ref{Thm:Kernel}. Specifically, we use quadrature along the contours $\Gamma_j$ given in Appendix~\ref{a:regularity}.  Since the precise paths of steepest descent do not need to be followed, we use piecewise-affine contours such that the angle of the contour that passes through each $\kappa_j$ agrees with the local path of steepest descent.  The routines in \cite{RHPackage} provide a robust framework for visualizing and computing such contour integrals.  In general, Clenshaw--Curtis quadrature is used on each affine component. To ensure accuracy for arbitrarily large $s$, the contour that passes through $\kappa_j$ is chosen to be of length proportional to $1/\sqrt{|s \omega_t''(\kappa_j)|}$.  This ensures that the Gaussian behavior near the stationary point is captured accurately in the large $s$ limit.  If all 
deformations are performed correctly, with this scaling behavior, a fixed number of sample points for Clenshaw--Curtis quadrature can be used for all $s$.  
A more in-depth discussion of this idea is given in \cite{TrogdonUnifiedNumerics} and \cite{TrogdonThesis}.

For reference purposes, the above method should be compared to a more restricted approach for the computation of generalized Airy functions 
presented in \cite{chinhedstrom}.  The authors of this paper compute special functions which correspond to $\omega(k) = k^p/p - i k^q/q$ for $m =-1,0$, i.e., they introduce dissipation into their special functions which corresponds to adding artificial viscosity into a finite-difference scheme for a hyperbolic system.  With this artificial dissipation they are able to characterize the asymptotic behavior of finite-difference schemes in terms of these special functions.

\paragraph{Example: Airy function.}

\begin{figure}[tbp]
\centerline{\includegraphics[width=.45\linewidth]{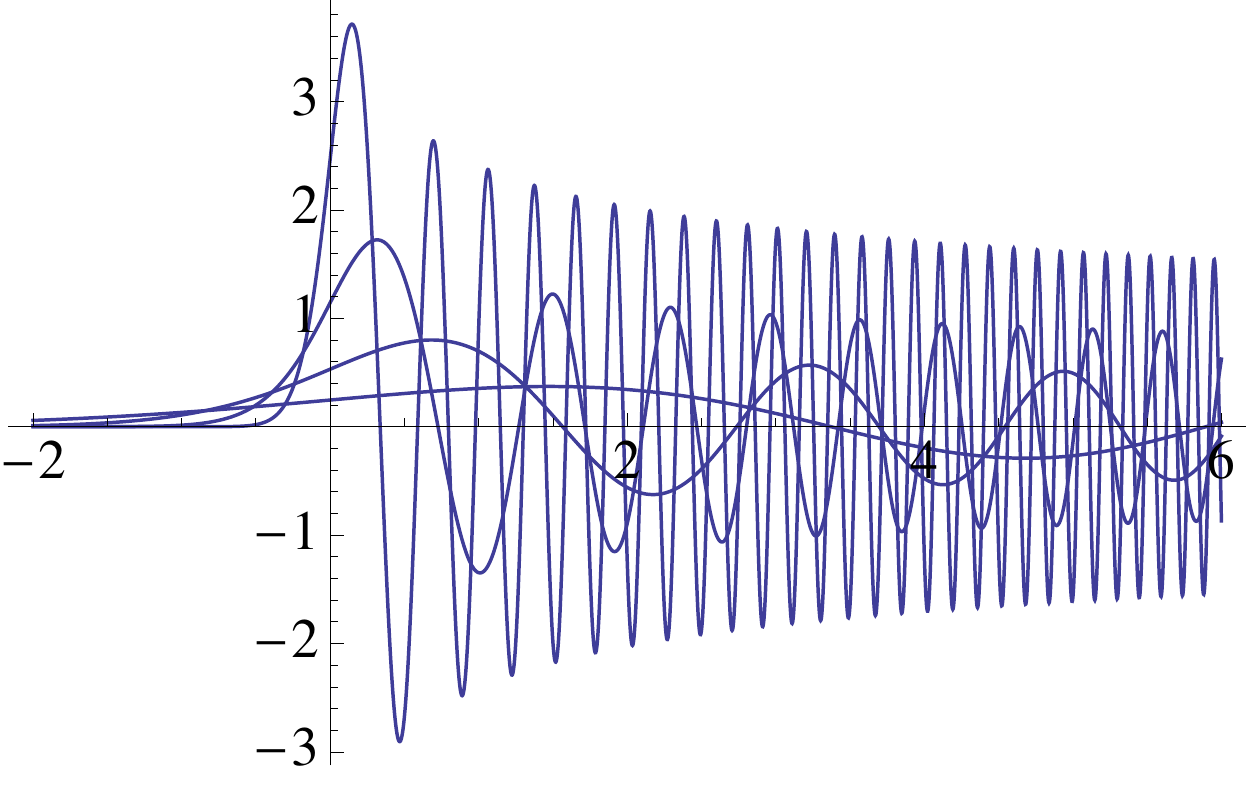}\includegraphics[width=.45\linewidth]{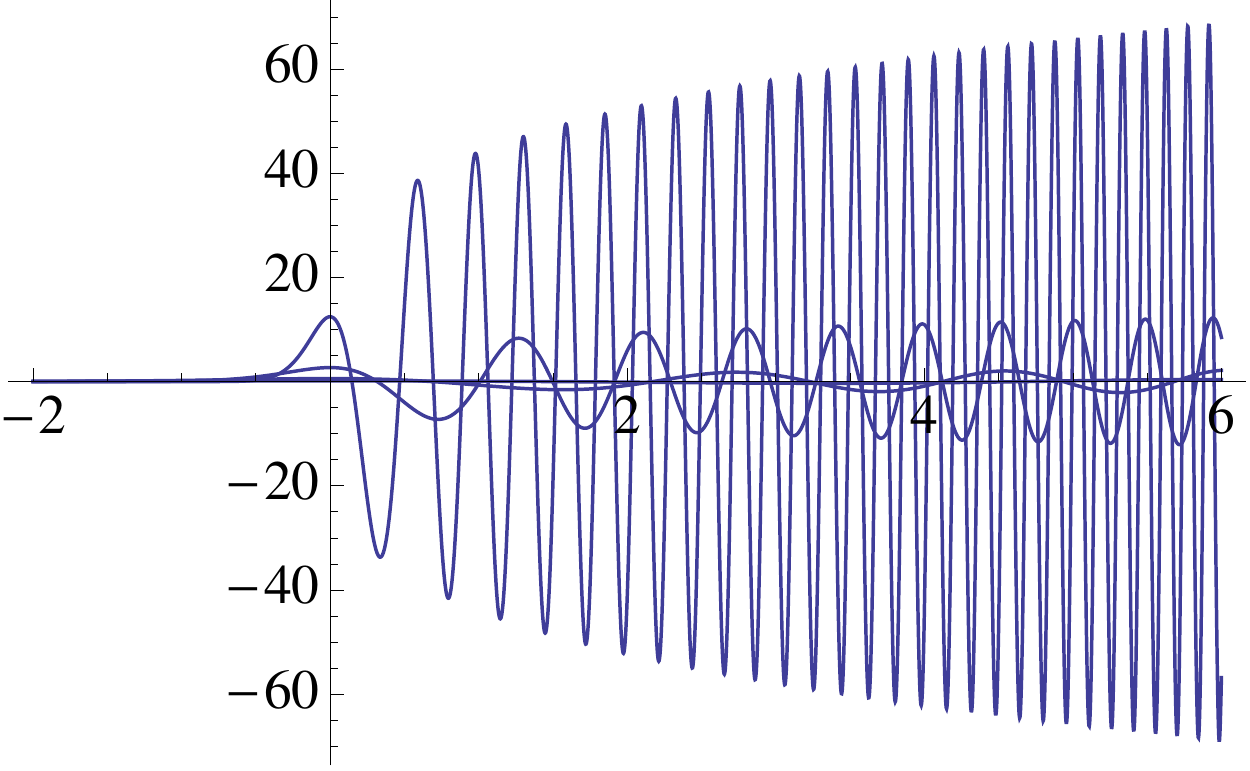}}
\caption{Plots of $I_{\omega,m}(y,t)$ with $\omega(k) = k^3$ versus $y$ for $t = 1,0.1,0.01,0.001$.  
Left: The scaled Airy function ($m = 0$). 
Right: The first derivative of the scaled Airy function ($m=1$).}
\label{Figure:Airy}
\end{figure}

When $\omega(k) = k^3$, the functions $I_{\omega,m}(y,t)$ are scaled derivatives and primitives of the Airy function. 
This function is displayed in Figure~\ref{Figure:Airy} for various values of $t$.   
See also Fig.~\ref{f:stokes},
where a primitive of the scaled Airy function ($m = -1$) was shown.
[But note that in Fig.~\ref{f:stokes} the dispersion relation was $\omega(k) = -k^3$, which results in a switch $y\mapsto -y$.]
It is clear that while the Airy function is bounded, its derivative grows in $x$.  
This is in agreement with Theorem~\ref{Thm:Kernel}.

\paragraph{Example: A higher-order solution.}

\begin{figure}[tbp]
\centerline{\includegraphics[width=.45\linewidth]{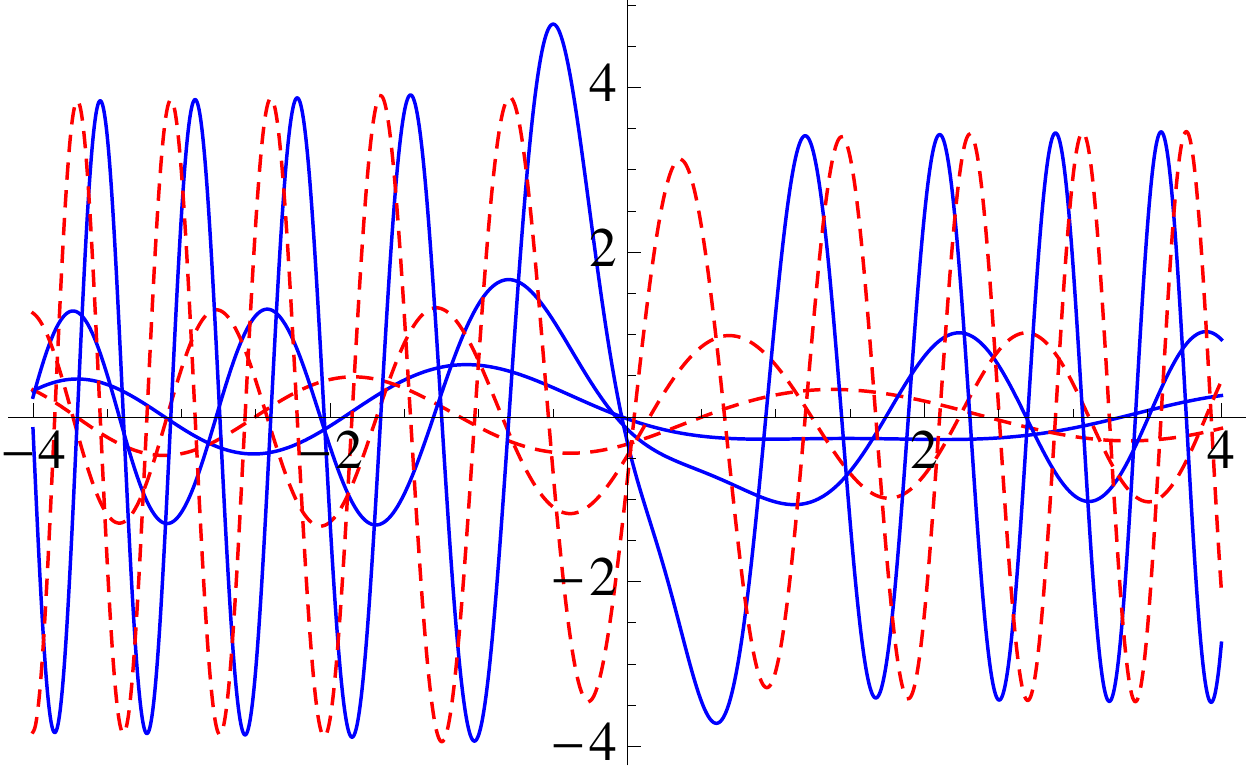}\includegraphics[width=.45\linewidth]{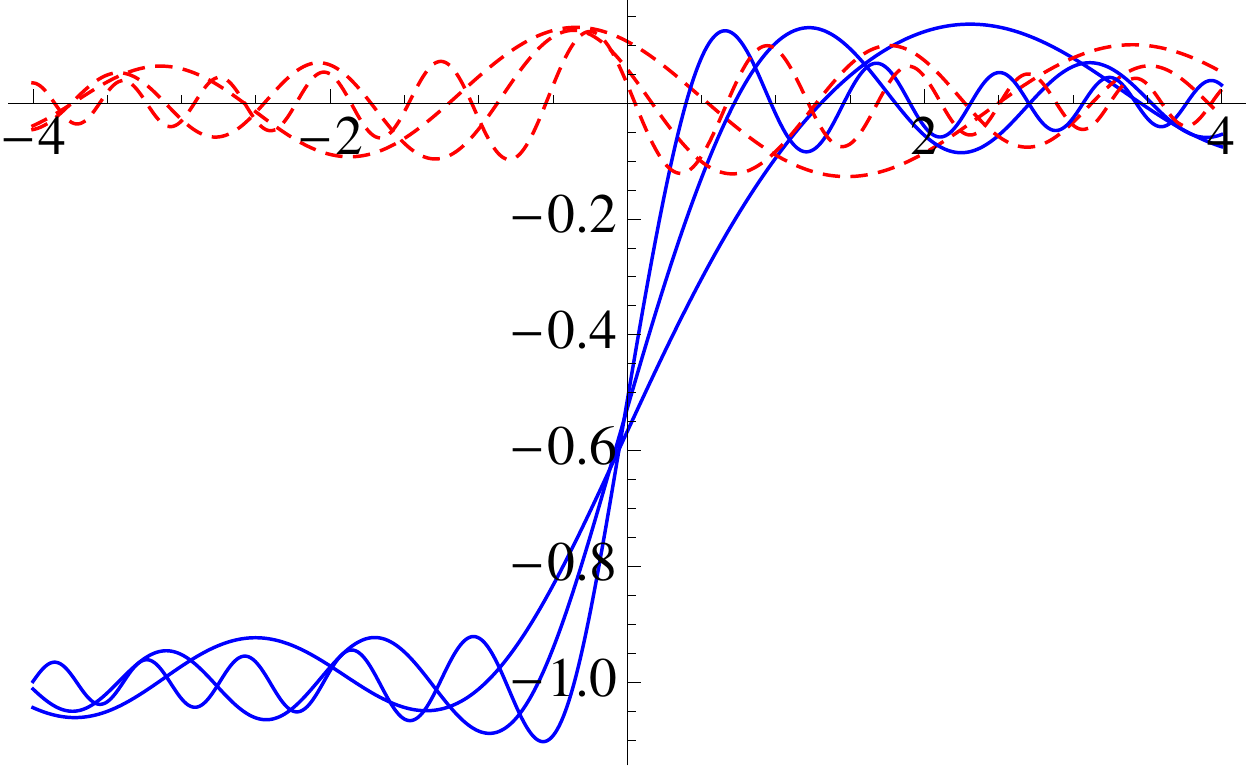}}
\caption{Plots of $I_{\omega,m}(y,t)$ with $\omega(k) = k^4 + 2 k^3$ versus $y$ for $t = 0.1,0.01,0.001$ (solid: real part, dashed: imaginary part).
Left: $m = 1$.  Right: $m = -1$.}
\label{Figure:Airy2}
\end{figure}

When the dispersion relation is non-monomial, the situation is more complicated. 
Consider for example $\omega(k) = k^4 + 2k^3$.
In this case $I_{\omega,m}(y,t)$ is no longer a similarity solution.  
Furthermore, it has non-zero real and imaginary parts.   
This function is displayed in Figure~\ref{Figure:Airy2} for various values of $t$.

\section{Concluding remarks}
\label{s:discussion}

We have obtained an asymptotic expansion for the short-time asymptotics of the solution of linear evolution PDEs with discontinuous ICs, including precise
error estimates.
The results apply to generic ICs (i.e., non-piecewise constant, non-compact support).
Moreover, the results extend to
arbitrary dispersion relations, 
multiple discontinuities,
and discontinuous derivatives of the IC.
In a forthcoming publication we will show that these results are also instrumental to characterize 
discontinuous BCs and corner singularities in IBVPs using the unified approach presented in \cite{Fokas2008}.
We end this work with a further discussion of the results.

{\bf1.} 
We have shown that the short-time asymptotic behavior of the solution of an evolution PDE with singular ICs is governed by similarity solutions
and classical special functions.
This is analogous to what happens in the long-time asymptotic behavior.
In that case, however, it is the discontinuities of the Fourier transform that provide the singular points for the analysis
(in addition of course to the stationary points or saddle points characteristic of the PDE).
In turn, these are related to the slow decay of the ICs at infinity.
In this sense, the short-time and long-time behavior are dual expressions of the characteristic behavior of a linear PDE.

{\bf2.} 
We have also shown that the solutions of dispersive linear PDEs exhibits Gibbs-like behavior in the short-time limit. This Gibbs-like behavior is robust, meaning that it persists under perturbation.
To explain this point, one should consider the obvious question of what happens with ICs which are a ``smoothed out'' discontinuity,
namely, a sharp but continuous transition from one value to a different one.
Such an IC can be considered to be a small perturbation of a step discontinuity in $L^2(\mathbb R) \cap L^1((1+|x|)^\ell)$.
Thus, as long as the IVP is well-posed,
the continuous dependence of the solution of the IVP on the ICs implies that 
a small change in the ICs will only produce a small change in the solution.

Let us briefly elaborate on this point.
Obviously if the perturbed IC is continuous, the solution of the PDE will converge uniformly to it as $t\downarrow0$.
Therefore, the Gibbs phenomenon that is present for the unperturbed solution will eventually disappear in the perturbed solution in this limit.
On the other hand, in Appendix~\ref{a:approx} we show that, 
if the perturbation is sufficiently small, one can still expect to observe a similar Gibbs-like effect \textit{at finite times}.

{\bf3.} 
The Gibbs-like behavior has been noticed in a couple of cases for nonlinear PDEs.
In particular, DiFranco and McLaughlin \cite{DiFranco2005} studied the behavior of the defocusing nonlinear Schr\"o\-dinger (NLS) equation with box-type IC.  The semiclassical focusing NLS equation was considered in \cite{jenkins} by Jenkins and McLaughlin.  Kotlyarov and Minakov \cite{KotlyarovMinakov} studied the behavior of the Korteweg-de Vries (KdV) and modified KdV equations with Heaviside ICs. 
In both cases, these authors showed that the behavior of the nonlinear PDE for short times is given to leading order by the behavior of the linear PDE.
And in both cases, in order to characterize the phenomenon it was necessary to use complete integrability of the nonlinear PDEs,
as well as Deift and Zhou's nonlinear analogue of the steepest descent method for oscillatory Riemann-Hilbert problems \cite{CPAM47p199,AM137p295}.
But the results of this work make it clear that this behavior: 
(i) is not a nonlinear phenomenon, and it also applies to linear PDEs;
(ii) is a general phenomenon, not limited to a few special PDEs.

{\bf4.}
At the same time it is true that for many nonlinear PDEs the nonlinear terms require $O(1)$ times in order to produce an
appreciable effect on the solution.
Therefore it is reasonable to expect that the results of this work will also provide the leading-order behavior of the solution
of many nonlinear PDEs for short times.
Indeed, Taylor \cite{Taylor2006} studied a generalized NLS equation (which is not completely integrable), 
and again characterized the behavior of the solutions for short times in terms of those of the linearized PDE.
It is hoped that such results can be generalized to other kinds of nonlinear PDEs.

{\bf5.} 
Of course, for larger times the solutions of linear and nonlinear PDEs with discontinuous ICs are very different from each other:
While for linear PDEs the oscillations spread out thanks to the similarity variable, 
for nonlinear PDEs the discontinuity gives rise to dispersive shock waves (DSWs);
namely, an expanding train of modulated elliptic oscillations with a fixed spatial period, whose envelope interpolates 
between the values of the solution at either side of the jump.
Such a nonlinear phenomenon has been known since the 1960's \cite{gurevichpitaevski},
and a large body of work as been devoted to its study (e.g., see \cite{biondinikodama2006,dubrovin2006,EgorovaTeschl,elkrylov1995,gravaklein2007,hoefer2006,kamchatnov,kamvissismclaughlinmiller,kodama1999} and references therein).
To the best of our knowledge, however, such behavior was never compared to the corresponding one for linear PDEs,
unlike what was done for the long-time asymptotics (e.g., see \cite{AS1981,NMPZ1984}).

{\bf6.} 
We reiterate that this Gibbs-like behavior of dispersive PDEs
\textit{is not a numerical artifact of a numerical approximation to the solution of the PDE,
but it instead a genuine feature of the solution itself}. 
We believe that this is perhaps the most important result of this work, 
since it has concrete implications for numerical analysis and the numerical solution of dispersive IVPs.
Namely, 
when performing numerical simulations of dispersive PDEs, 
one must be very careful to distinguish among spurious Gibbs features induced by the truncation of a Fourier series representation,
spurious Gibbs oscillations generated by numerical dispersion (introduced by the numerical scheme used to solve the PDE),
and actual Gibbs-like behavior generated by the PDE itself.



{\bf7.} 
From a philosophical point of view, one may ask why consider PDEs with discontinuous ICs at all.  
In this respect we note on one hand that, apart from any physical considerations,
studying these kinds of ICs is important from a mathematical point of view to understand the properties of the PDE and its solutions. 
Also, on the other hand, such a study also makes perfect sense physically.  
For example, one only need think about hyperbolic systems, 
for which considerable effort is devoted to the study of shock propagation.
These shocks are discontinuities in the solution,
and describe actual physical behavior.
Even though such discontinuities are only approximation of a thin boundary layer,
the fact remains nonetheless that representing such situations with discontinuous solutions
is a convenient mathematical representation of the actual physical behavior.
More in general, while the PDE holds in the interior of the domain $(x,t)\in\Real\times\Real^+$,
the IC is posed on the \textit{boundary} of this domain.
In this sense $t=0$ is always a singular limit. 
Indeed, the results of Section~\ref{s:regularity} show that, generally speaking, 
the solution on the interior of the domain is smooth
even when the IC is singular.

        

\appendix

\section{Brief review of well-posedness results}
\label{a:wellposedness}

In this section we briefly review some well-known results about well-posedness of the IVP for the PDE~\eref{e:pde} with 
dispersion relation~\eref{e:dispersion} and IC~\eref{e:IC}. We define the Fourier transform pair for $f \in L^2(\mathbb R)$ by
\begin{align}\label{e:ft-pair}
\hat f(k) = \int_{-\infty}^\infty e^{-ikx} \hat f(x)\d x, \quad f(x) = \frac{1}{2\pi} \int_{-\infty}^\infty e^{ikx} \hat f(k)\d x.
\end{align}
Throughout, we use the caret ($~\hat {}~$) to denote the spatial Fourier transform.

\begin{definition}\label{d:classical}
The function $q(x,t)$ is a classical solution of the PDE \eref{e:pde} with~\eref{e:dispersion}
in an open region $\Omega\subset \Real^2$ if all derivatives present in the PDE exist for all $(x,t) \in \Omega$ and the PDE is satisfied pointwise.
\end{definition}

Recall that $\omega(k)$ needs to satisfy certain conditions in order for the IVP for~\eref{e:pde} to be well posed. Specifically,
it is straightforward to see that $\Im[\omega(k)]$ must be bounded from above.  Letting $n = \deg[\omega(k)]$ this condition implies
$\Im\omega_n\le0$ if $n$ is even and $\omega_n\in\Real$ if $n$ is odd.   Also recall that the PDE is said to be dispersive if $\omega''(k)\ne0$ \cite{whitham}.

\begin{definition}
A function $q(x,t)$  is a weak solution of~\eref{e:pde} with~\eref{e:dispersion} in an open region $\Omega$ if
\begin{align}\label{e:weak-form}
L_\omega[q,\phi] = \int_\Omega q(x,t)(-i \partial_t\phi(x,t) - \omega(i \partial_x) \phi(x,t)) \d x \d t = 0,
\end{align}
for all $\phi \in C^\infty_c(\Omega)$ (with the subscript $c$ denoting compact support).
\end{definition}

\begin{definition}
A function $q(x,t)$ is an $L^2$ solution of the IVP for \eref{e:pde} with dispersion relation~\eref{e:dispersion} and IC~\eref{e:IC} if:
(i) $q \in C^0([0,T];L^2(\mathbb R))$, 
(ii)~$q$~satisfies \eqref{e:weak-form} with $\Omega = \mathbb R \times \mathbb R^+$, and 
(iii) $q(\cdot,0) = q_o$ a.e.
\end{definition}

We now show that the function $q(x,t)$ defined by the Fourier transform reconstruction formula
\eref{e:soln-form}
with
$\theta(x,t,k) = kx-\omega(k)t$
is an $L^2$ solution of the IVP provided the imaginary part of $\omega(k)$ is bounded above and $q_o \in L^2(\mathbb R)$.  
To see this, one can use the convolution property of the Fourier transform, which is a consequence of the Plancherel theorem: 
if $f, g \in L^2(\mathbb R)$, then
\begin{align}\label{e:planch}
\int_{\mathbb R} f(x) g(x) \d x = \frac{1}{2\pi} \int_{\mathbb R} \hat f(k) \hat g(-k) \d k.
\end{align}
Applying~\eqref{e:planch} (in $x$) to \eqref{e:weak-form} yields
\[
L_\omega[q,\phi]= \frac{1}{2\pi}\int_{\mathbb R^+}\int_{\mathbb R} e^{-i\omega(k)t} \hat q_o(k) (-i \partial_t \hat \phi(-k,t) - \omega(k) \hat \phi(-k,t)) \d k \d t.
\]
But note that 
\begin{align*}
-i e^{-i\omega(k)t} \hat q_o(k)  \partial_t \hat \phi(-k,t) = -i \partial_t(e^{-i\omega(k)t} \hat q_o(k) \hat \phi(-k,t)) + \omega(k) e^{-i\omega(k)t} \hat q_o(k) \hat \phi(-k,t)\,.
\end{align*}
From this it follows that $L_\omega[q,\phi] = 0$ because of the compact support of~$\phi$.

We next show that the Fourier transform solution is unique.
To see this, take $\phi(x,t) = X(x)T(t)$, for any $L^2$ solution we have
\begin{align*}
L_\omega[q,\phi]= \frac{1}{2\pi}\int_{\mathbb R}&\hat X(-k) \int_{\mathbb R^+} \hat q(k,t) (-i \partial_t T(t) - \omega(k) T(t)) \d t \d k = 0.
\end{align*}
Since the inner integral defines a locally integrable function (with polynomial growth in $k$) and $X$ is arbitrary, 
the inner integral must vanish for a.e.~$k$.  
Specifically, this follows from the density of $C_c^\infty(\mathbb R)$ in the Schwartz class $\mathcal S(\mathbb R)$. 
This is the weak form of an ODE for $\hat q(k,t)$ and we must show that the obvious solution of this is the only solution. We rewrite this condition as (for fixed $k$)
\begin{align}\label{e:t-int}
0 = -i \int_{\mathbb R^+} e^{i\omega(k)t} \hat q(k,t) (e^{-i\omega(k)t}  T(t))_t \d t.
\end{align}
Now, \eqref{e:t-int} implies that the integral of $e^{i\omega(k)t} \hat q(k,t)$ against any $C_c^\infty(\mathbb R^+)$ function with integral zero is zero: a $C_c^\infty(\mathbb R^+)$ function has integral zero if and only if it is the derivative of $C_c^\infty(\mathbb R^+)$ function.   Now let $\phi, \psi \in C_c^\infty(\mathbb R^+)$ and choose $\psi$ so that it integrates to one.  Then
\begin{align*}
\phi(t) - \psi(t) \int_{\mathbb R^+} \phi(s)ds,
\end{align*}
is a test function that integrates to zero.  We find
\begin{align*}
\int_{\mathbb R^+} e^{i\omega(k)t} \hat q(k,t) \phi(t) \d t = \int_{\mathbb R^+} \phi(s)ds \cdot \int_{\mathbb R^+} e^{i\omega(k)t} \hat q(k,t) \psi(t) \d t,
\end{align*}
and the inner integral in the right-hand side 
must be a constant $c(k)$ (independent of $\psi$). 
Thus
\begin{align*}
\int_{\mathbb R^+} [e^{i\omega(k)t} \hat q(k,t)-c(k)] \phi(t) \d t = 0,
\end{align*} 
for all $\phi \in C_c^\infty(\mathbb R^+)$ and $e^{i\omega(k)t} \hat q(k,t) = c(k)$ is constant for a.e. $t$. This proves $\hat q(k,t) = e^{-i\omega(k)t} \hat q_o(k)$.
Finally, examining \eqref{e:soln-form}, it is easily seen that 
the solution, as a function in $C^0([0,T];L^2(\mathbb R))$ depends continuously on the initial data.

\section{Asymptotics of the special functions and regularity results}
\label{a:regularity}

In the first part of this section we concentrate on the steepest descent analysis of $I_{\omega,m}(x,t)$, 
which will also give us results concerning an important convolution kernel $K_t(x) \triangleq I_{\omega,-1}(x,t)$.  
Then in the last part of the section we apply these results to the IVP for the linear PDE~\eref{e:pde}.

We are interested in the asymptotics of $I_{\omega,m}(x,t)$ in two particular limits.  
First, in this section we need estimates for fixed $t$ as $|x| \rightarrow \infty$.  Estimates for fixed $x$ as $t \rightarrow 0^+$ are also needed.  For $|x| > 0$, $t > 0$, we rescale, by setting $\sigma = \sign(x)$, $k = \sigma (|x|/t)^{1/(n-1)}z$
\begin{align}\label{e:rescale}
\begin{split}
I_{\omega,m}(x,t) &= \frac{1}{2 \pi}\sigma^m \left(\frac{|x|}{t}\right)^{-m/(n-1)} \int_{C} e^{X(iz-i\omega_n \sigma^n z^n-i R_{|x|/t}(z))} \frac{\d z}{(iz)^{m+1}},\\
R_{|x|/t}(z) &= \sum_{j=2}^{n-1} \omega_j \left(\frac{|x|}{t}\right)^{\frac{j-n}{n-1}} (\sigma z)^j, \quad X = |x| \left( \frac{|x|}{t} \right)^{1/(n-1)}.
\end{split}
\end{align}
Here analyticity allows the deformation back to $C$ (see Figure~\ref{f:contour}) after the change of variables.

The benefit of this scaling, as we will see, is that $R_{|x|/t}$ has coefficients that decay as $|x|/t$ increases. For $x < 0$ we perform a negative scaling to get it in a form where $X > 0$.  This effectively maps $\omega_n$ to $ \sigma^n \omega_n$ and because the lower-order terms in $R_{x/t}$ create lower-order effects the computation proceeds in the same way.  Therefore, we provide all results for $x > 0$ noting the correspondence to $x < 0$. Define
\begin{align*}
\Phi_{|x|/t}(z) = iz-i\omega_n \sigma^n z^n-i R_{|x|/t}(z),
\end{align*}
where $\{z_j\}_{j=1}^{n-1}$ are the roots of $\Phi_{|x|/t}'(z) = 0$ ordered counter-clockwise from the real axis.  It is clear that there exists $n-1$ distinct roots for sufficiently large $|x|/t$ and such an ordering is possible.

Because it is the liming case, consider $R_{|x|/t} \equiv 0$ (\emph{i.e.} a monomial dispersion relation).  We are interested in the roots of
\begin{align*}
 n\omega_n \sigma^n z^{n-1} = 1.
\end{align*}
\begin{itemize}
\item  If $n$ is even we have one root on the real axis and $(n-2)/2$ roots in the upper-half plane.
\item If $n$ is odd and $\omega_n\sigma^n$ is positive we have two roots on the real axis and $(n-3)/2$ roots in the upper-half plane.
\item If $n$ is odd and $\omega_n\sigma^n$ is negative we have no roots on the real axis and $(n-1)/2$ roots in the upper-half plane.
\end{itemize}
It is straightforward to check that for sufficiently large $|x|/t$ these statements hold for $\Phi_{|x|/t}'(z) = 0$.  Define $N(n)$ to be this number of roots in the closed upper-half plane.

Consider the region $D = \{ k : \Im  \omega_n \sigma^n k^n > 0\}$.  This is the region in which $e^{-i\omega_n \sigma^n k^n}$ is unbounded and any contour deformation should avoid this region for large $k$.  It is straightforward to check that $D$ consists of $n$ wedge-like sectors emanating from the origin.  The steepest descent path though $z_j$ satisfies
\begin{align*}
0 = \Im \Psi_{|x|/t,j}(z),  \quad \Psi_{|x|/t,j}(z) = \Phi_{|x|/t}(z) - \Phi_{|x|/t}(z_j).
\end{align*}
Writing $z = r e^{i \theta(r)}$ for large $r$ we find
\begin{align*}
\cos( n \theta(r)) + \bigo(r^{1-n} + r^{-1} (|x|/t)^{-1/n})) = 0.
\end{align*}
Therefore using analyticity of the inverse cosine function near a zero of cosine
\begin{align*}
\theta(r) = \frac{2m + 1}{2n} \pi + \bigo(r^{-1} ( 1 + (|x|/t)^{-1/n})).
\end{align*}
We note that the steepest descent directions of $e^{-i \omega_n \sigma^n k^n}$ are given by a subset of $\theta = \frac{2m+1}{2n} \pi$, $m = 0, \ldots, 2n-1$ such that $\omega_n \sigma^n \sin n \theta > 0$.  Thus, in this sense any unbounded portions of steepest descent paths are asymptotic, uniformly in $|x|/t$, to a steepest descent path of $e^{-i \omega_n \sigma^n k^n}$.  We work towards understanding the paths of steepest descent that pass through $\{z_j\}_{j=1}^{N(n)}$.    Next, we note that in the monomial case ($R_{|x|/t} \equiv 0$)
\begin{align*}
\Im (iz_j - i \omega_n \sigma^n z_j^n) =  \frac{n-1}{n} \Re z_j.
\end{align*}
Thus there can be path of steepest descent or ascent that connects to stationary points only if they have equal real parts.  So, we compare their real parts:
\begin{align}\label{e:realpart}
\Re (iz_j - i \omega_n \sigma^n z_j^n) = - \frac{n-1}{n} \Im z_j.
\end{align}
It is clear that the exponent evaluated at stationary points in the upper-half plane has a smaller real part.  Thus, any stationary point in the upper-half plane has no steepest descent path that connects to any other stationary point. Finally, it can be shown that the steepest descent path through $z_j$ must be asymptotic to the closest (with respect to argument) steepest descent paths of $e^{-i\omega_n \sigma^n k^n}$.  We establish the following:
\begin{lemma}\label{Lemma:Gammaj}
For sufficiently large $|x|/t$ there exists unique, disjoint contours $\Gamma_j \subset \{z : \Im \Phi_{|x|/t}(z) = \Im \Phi_{|x|/t}(z_j)\}$, $j = 1, \ldots, N(n)$ such that 
\begin{itemize}
\item $z_j \in \Gamma_j$,
\item $\Gamma_j$ corresponds to the path of steepest descent from $z_j$, and
\item $\Gamma_j$ is asymptotic in each direction to a steepest descent path of $e^{-i\omega_n \sigma^n k^n}$.
\end{itemize}
\begin{proof}
Previous arguments demonstrate this in the monomial case.  For the general case we note that $\{z_j\}$ converge to roots of $1 = n \omega_n \sigma^n z^n$ as $|x|/t \rightarrow \infty$ and the same conclusions follow.
\end{proof}
\end{lemma}

From \eref{e:realpart}, in the monomial case, $z_j \in \overline{D}$.  Furthermore, $z_j$ and $z_{j+1}$ lie in distinct sectors $D_j$ and $D_{j+1}$ of $D$ and one and only one sector $H_j$ of $D^c$ lies between these two sectors (with counter-clockwise ordering).  Let $H_0$ be the sector of $D^c$ that lies before $z_1$ and $H_{N(n)}$ be the sector of $D^c$ that lies after $z_{N(n)}$ with the same ordering.  It is also clear that  $\Gamma_j$ cannot limit to infinity in any sector besides $H_{j-1}$ and $H_{j}$ as this would imply that the imaginary part of the exponent varied along the path. Next we understand the change of variables that is used along the steepest descent path.



It is important that $\Gamma_j$ passes through one and only one stationary point $z_j$ of the new exponent.  Further, in the case that $R_{|x|/t} \equiv 0$ (\emph{i.e.} a monomial dispersion relation) it is clear how to proceed: a straightforward application of the method of steepest descent for integrals will give the leading-order term.  A derivation of the result in the general case requires some technical work.  Define the variable $v_{|x|/t}(s)$ by the equation
\begin{align*}
\frac{\Psi_{|x|/t,j}(z_j + s v_{|x|/t})}{s^2} + 1 = 0, ~~ { v_{|x|/t}(0)= \pm(-\half \Psi_{|x|/t,j}''(z_j))^{-1/2},}
\end{align*}
{ where $\pm$ is chosen so that $\real \Psi_{|x|/t,j}(z_j + s v_{|x|/t}) \leq 0$.}   We use $v_\infty$ to refer to the case where $R_{|x|/t}\equiv 0$. The Implicit Function Theorem can be applied for each $s \in \mathbb R$, producing the function $v_{|x|/t}(s)$ which depends smoothly on $s$ and $|x|/t$ provided $|x|/t$ is sufficiently large, so that the coefficients of $R_{|x|/t}$ are sufficiently small.  Applying the change of variables $k = \tau_{|x|/t}(s) = z_j + s v_{|x|/t}(s)$ we find
\begin{align*}
\int_{\Gamma_j} e^{X(iz-i\omega_n \sigma^n z^n-i R_{|x|/t}(z))} \frac{\d z}{(iz)^{m+1}} = \int_{\mathbb R} e^{-Xs^2} \frac{\tau'_{|x|/t}(s)}{(i\tau_{|x|/t}(s))^{m+1}} ds.
\end{align*}
Now, one should expect that a Taylor expansion of
\begin{align*}
F_{m,|x|/t}(s) = \frac{\tau'_{|x|/t}(s)}{(i\tau_{|x|/t}(s))^{m+1}} 
\end{align*}
at $s = 0$ would produce a series expansion for the integral.  But, because of the $|x|/t$-dependence, extra work is required.

\begin{lemma}\label{Lemma:Constants}
There exists positive constants $C_{\ell,R}$ and $\epsilon_{\ell,R}$ (depending only on $\ell$, $R$ and  $\omega$) such that for $|t/x| < \epsilon_{\ell,R}$ and $|s| < R$
\begin{align*}
\sup_{s \in \mathbb R} |F_{m,|x|/t}^{(\ell)}(s)| < C_{\ell,R}.
\end{align*}
\begin{proof}
We begin by examining $v_{|x|/t}$ and its derivatives.  First, because of the Implicit Function Theorem, $v_{|x|/t}(s)$ and all its derivatives depend continuously on $s$ and $|x|/t$ for $|x|/t$ sufficiently large.  Furthermore, all derivatives limit to $v_\infty$ pointwise as $|x|/t \rightarrow \infty$.   For every $s \in \mathbb R$ there exists $\epsilon_s$ such that for $|s'-s| < \epsilon_s$ and $|t/x| < \epsilon_s$,
\begin{align*}
|v^{(\ell)}_{|x|/t}(s') - v^{(\ell)}_{\infty}(s')|  \leq |v^{(\ell)}_{|x|/t}(s') - v^{(\ell)}_{\infty}(s)| + |v^{(\ell)}_{\infty}(s) - v^{(\ell)}_{\infty}(s')| < \delta. 
\end{align*}
For $|s| \leq R$ we use compactness to cover $[-R,R]$ with a finite number of the balls $\{|s'-s_i| < \epsilon_{s_i}\}$, and let $\epsilon_{\ell,R} = \min_i \epsilon_{s_i}$.  It follows that for $|t/x|< \epsilon_{\ell,R}$, and $|s| \leq R$ that $|v^{(\ell)}_{|x|/t}(s) - v^{(\ell)}_{\infty}(s)| < \delta$. 
\end{proof}
\end{lemma}

\begin{figure}[tbp]
\centering
\includegraphics[width=.5\linewidth]{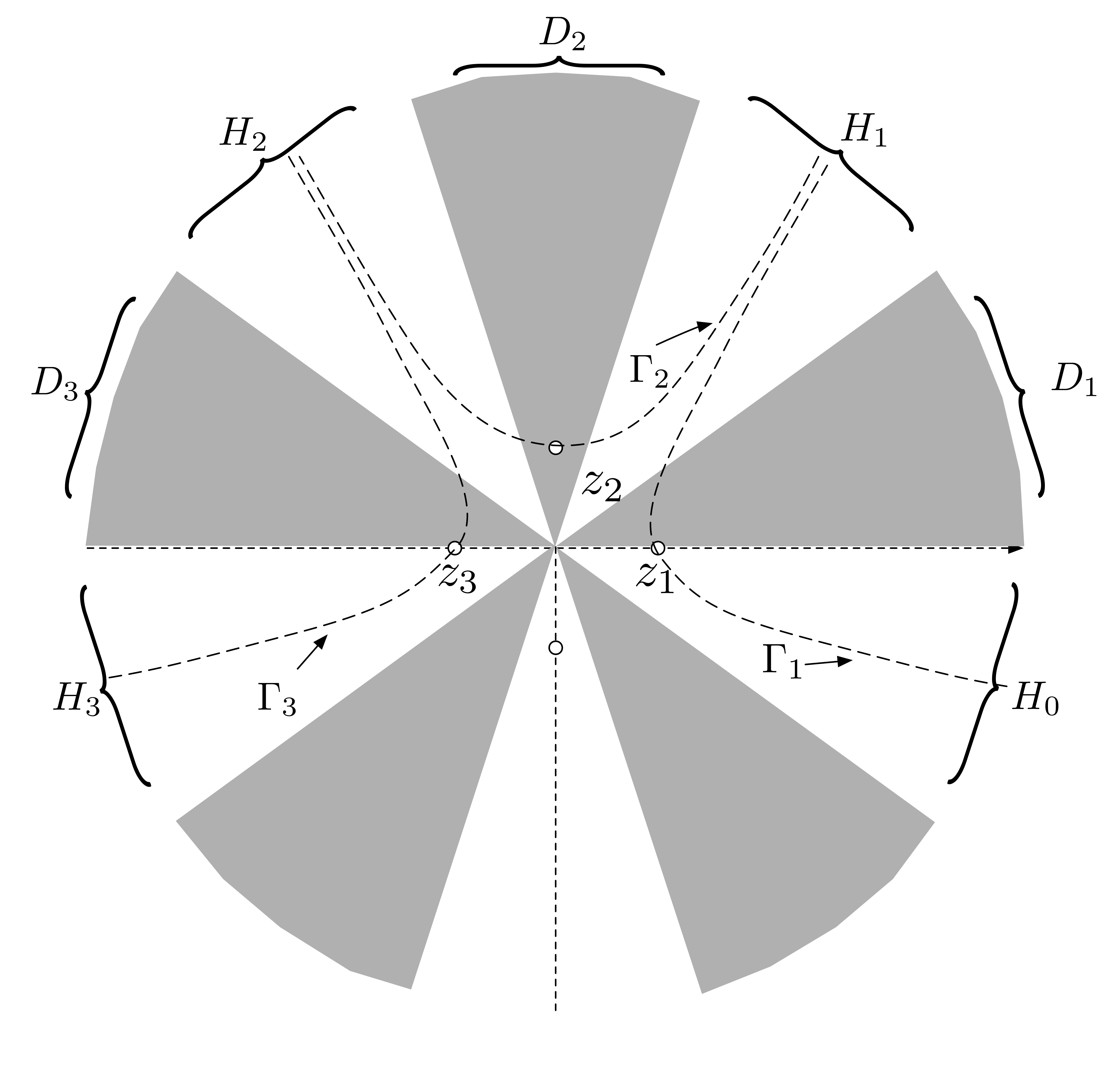}
\caption{\label{Figure:WhichGammaj} A schematic for $\omega(k) = k^5 + \bigo(k^{4})$.  The shaded region $S^c = \{z: \Im \omega_n \sigma^n z^n > 0 \}$ where the function $e^{-i\omega_n \sigma^n z^n}$ has growth. The circles represent the stationary points where $\Phi'_{|x|/t}(z) = 0$. The contours $\Gamma_j$ which are along the global paths of steepest descent.  This figure shows the definitions of the sectors $H_j$ and $D_j$.}
\end{figure}

From previous considerations and the convergence of $v_{|x|/t}$ to $v_{\infty}$ we have the following which is illustrated in Figure~\ref{Figure:WhichGammaj}.

\begin{lemma}\label{Lemma:WhichGammaj}
For sufficiently large $|x|/t$ and each $j = 1, \ldots, N(n)$, $z_j$ lies in a distinct sector $D_j \subset D$ and $\Gamma_j$ tends to infinity in both $H_{j-1}$ and $H_j$.  Furthermore, $H_{j}$ for $j = 2,\ldots, N(n)-1$ contains unbounded components of two contours $\Gamma_j$ and $\Gamma_{j+1}$, $H_0$ contains an unbounded component of $\Gamma_1$ and $H_{N(n)}$ contains an unbounded component of $\Gamma_{N(n)}$. 
\end{lemma}

This completes our characterization of the steepest descent paths and we consider the deformation of the integral.

\begin{figure}[tbp]
\centering
\subfigure[]{\includegraphics[width=.49\linewidth]{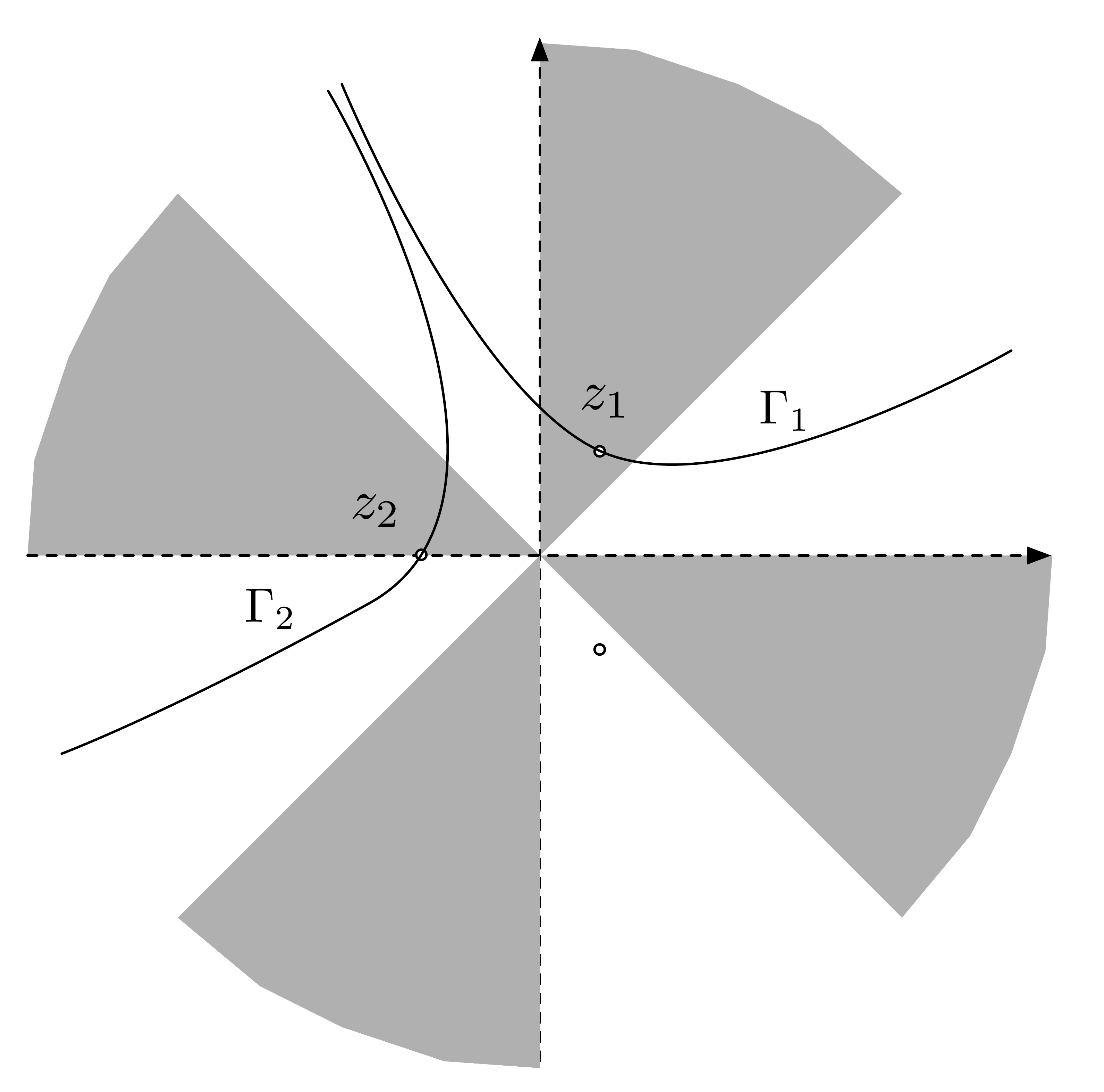}\label{Figure:Gammaj}}
\subfigure[]{\includegraphics[width=.49\linewidth]{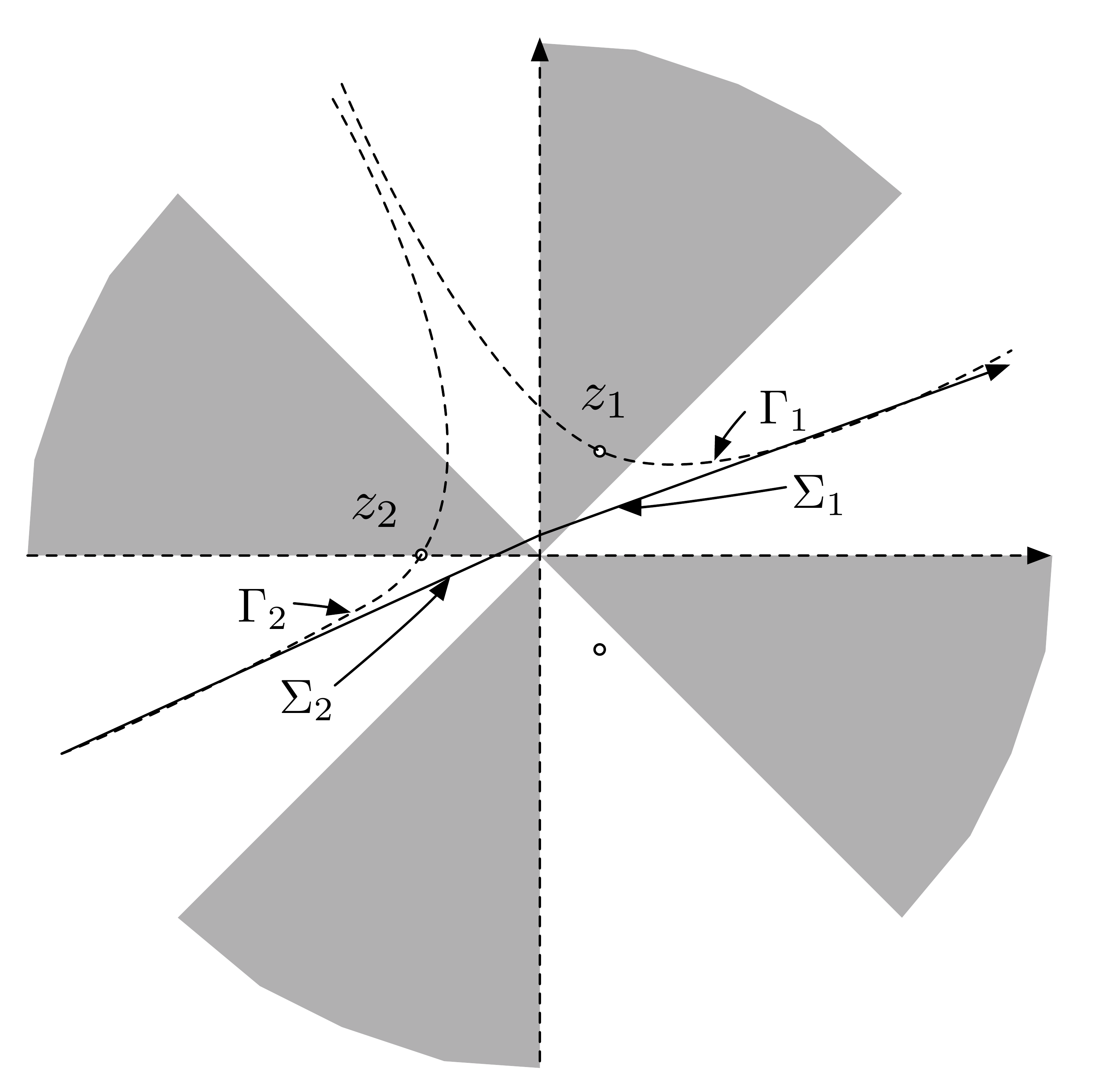}\label{Figure:Sigma}}
\caption{A schematic for $\omega(k) = -k^4 + \bigo(k^{3})$.  
The shaded region $S^c = \{z: \Im \omega_n \sigma^n z^n > 0 \}$ where the function $e^{-i\omega_n \sigma^n z^n}$ has growth. 
The circles represent the stationary points where $\Phi'_{|x|/t}(z) = 0$.  
(a) The contours $\Gamma_j$ which are along the global paths of steepest descent.  
(b) The initial deformation of the integral representation of $K_t(x)$, after scaling, to the contours $\Sigma_1$ and $\Sigma_2$.}
\bigskip
\centering
\subfigure[]{\includegraphics[width=.49\linewidth]{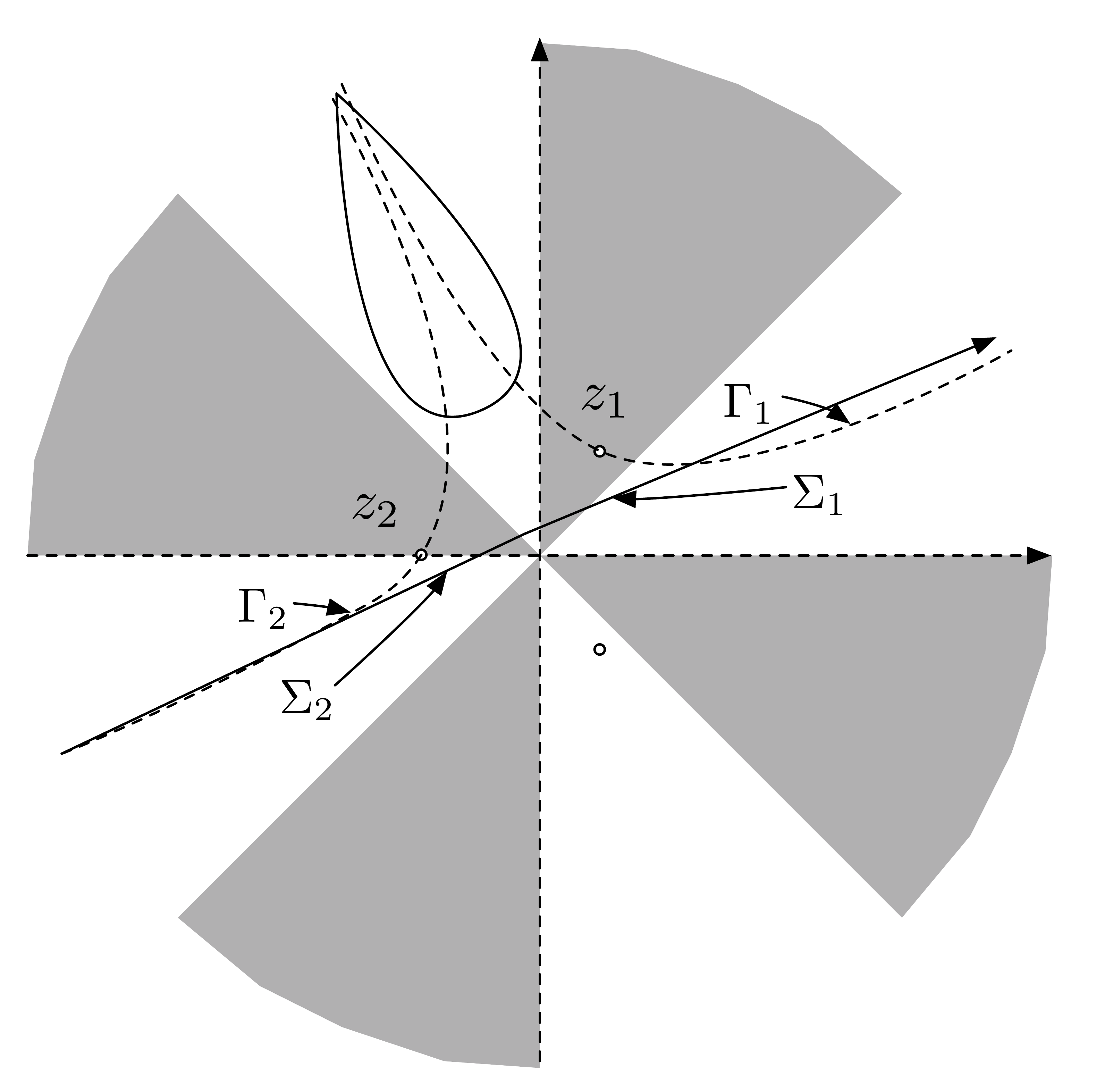}\label{Figure:Add}}
\subfigure[]{\includegraphics[width=.49\linewidth]{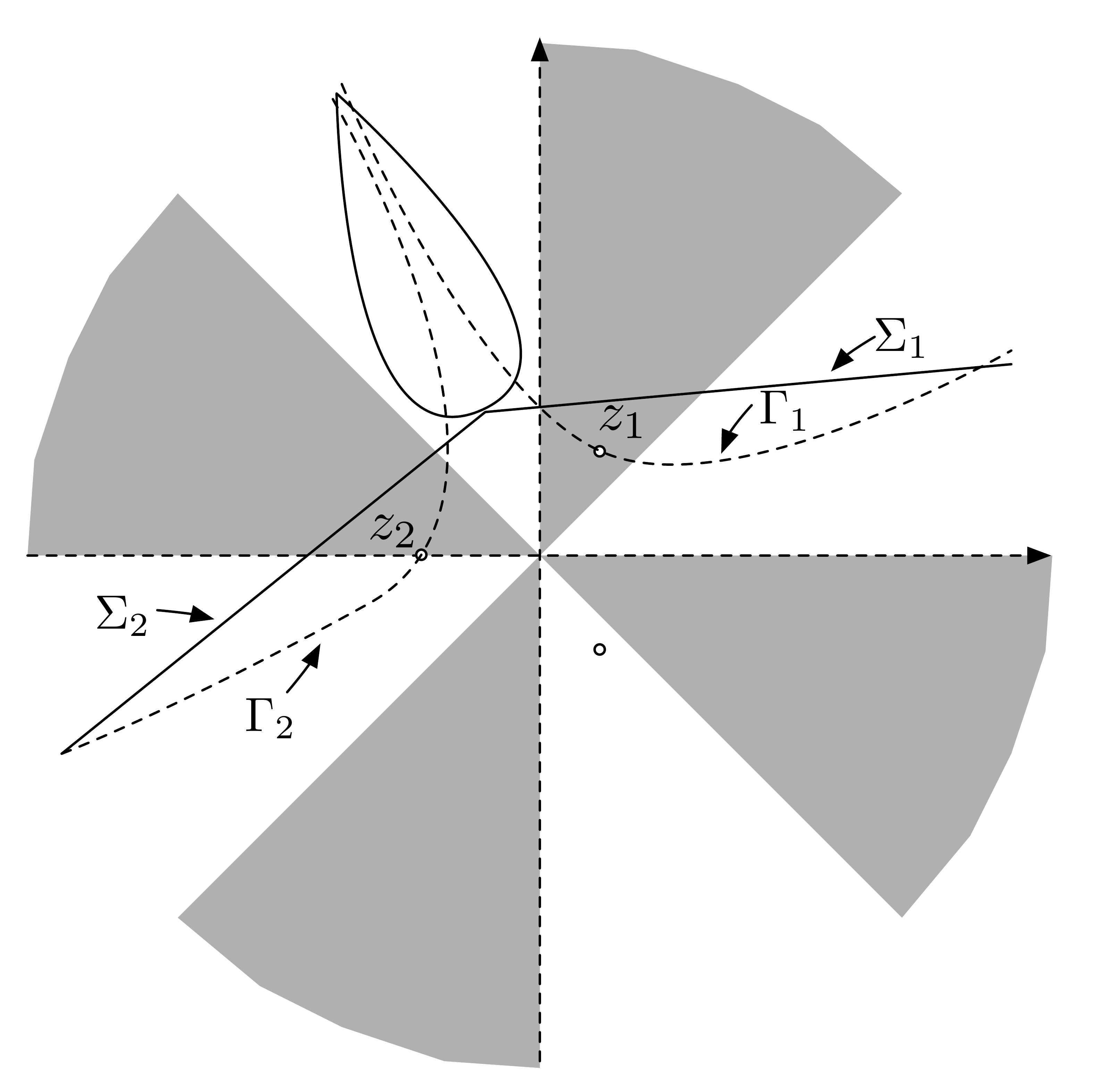}\label{Figure:Join}}
\caption{A a continuation schematic for $\omega(k) = -k^4$.  
The shaded region $S^c = \{z: \Im \omega_n \sigma^n z^n > 0 \}$ where the function $e^{-i\omega_n \sigma^n z^n}$ has growth. 
The circles represent the stationary points where $\Phi_{|x|/t}'(z) = 0$.  
(a) The addition of loop contours that contribute nothing (by Cauchy's Theorem) to the integral representation for $K_t(x)$.  
(b) The deformation of $\Sigma_1$ and $\Sigma_2$ to join with the loop contours.  At this point it is clear that these contours may be deformed to $\Gamma_1 \cup \Gamma_2$.}
\label{Figure:AddContours}
\end{figure}

\begin{lemma}
\label{Lemma:Deform}
For $m \geq -1$
\begin{align*}
\int_{C} e^{X\Phi_{|x|/t}(z)} \frac{\d z}{(iz)^{m+1}} = \sum_{j=1}^{N(n)} e^{X\Phi_{|x|/t}(z)} \int_{\Gamma_j}\frac{\d z}{(iz)^{m+1}},
\end{align*}
where $C$ is replaced with $\mathbb R$ if $m = -1$. 
\begin{proof}
We first work at the deformation of
\begin{align*}
\int_{1}^\infty e^{X\Phi_{|x|/t}(z)} \frac{\d z}{(iz)^{m+1}}
\end{align*}
off the real axis.  It can be seen that the boundary of the sector $H_0$ contains the real axis. It follows from the fact that $\Phi_{|x|/t}(z)$ has purely imaginary coefficients there exists an interval $[c,\infty)$ that is a subset of the boundary of the region $S = \{z: \Re \Phi_{|x|/t}(z)  \leq 0\}$ which may be above or below the real axis.  Furthermore, the component of $S$ whose boundary contains $[c,\infty)$ contains $\Gamma_1 \cap \{|z| > L\}$ for sufficiently large $L$ and $c$  can be taken to be independent of $|x|/t$. We justify the deformation of this integral to a contour that extends from $1$ to one of the points in $\Gamma_1 \cap \{|k| = L\}$ and then follows $\Gamma_1$ for $|k| > L$.  Call this contour $\Sigma_1$, see Figure~\ref{Figure:Sigma}. To establish this it suffices to demonstrate that
\begin{align}\label{large-R}
\int_{C_R}& e^{X(iz-i\omega_n \sigma^m z^n-i R_{|x|/t}(z))} \frac{\d z}{(iz)^{m+1}},\\
C_R& = \{ z : |z| = R, ~~ 0 \leq \pm \arg z \leq 1/(2n\pi) + \bigo(R^{-1}) \},
\end{align}
tends to zero for large $R$.  The $+,~-$ sign is taken if the deformation occurs above, below the real axis.  Lemma~\ref{Lemma:Gammaj} demonstrates the asymptotic form of $C_R$.   
Since the integrand itself does not decay uniformly when $m =-1$ we perform integration by parts.  Let $\gamma_R = \Gamma_1 \cap \{|k| = R\}$.  Then
\begin{multline*}
\int_{C_R} e^{X\Phi_{|x|/t}(z)} \d z = \left.\frac{iX^{-1}}{n \omega^n \sigma^n z^{n-1} + R'_{|x|/t}(z)} e^{X\Phi_{|x|/t}(z)} \right|_R^{\gamma_R}  
\\
+ iX^{-1} \int_{C_R}\left( \frac{n(n-1) \omega^n \sigma^n z^{n-2} + R''_{|x|/t}(z)}{(n \omega^n \sigma^n z^{n-1} - R'_{|x|/t}(z))^2} - \frac{i}{n \omega^n \sigma^n z^{n-1} - R'_{|x|/t}(z)} \right) e^{X\Phi_{|x|/t}(z)} \d z.
\end{multline*}
The boundary terms here drop out in the large $R$ limit.  The non-exponential factor in the integrand decays at least like $1/z^2$ so that it suffices to show the exponential is bounded on $C_R$ for sufficiently large $R$. This follows from the fact that for sufficiently large $R$, $C_R \subset S$.  It is clear that the argument also holds for $m > -1$. Similar reasoning may be applied to 
\begin{align*}
\int_{0}^\infty e^{X\Phi_{|x|/t}(z)} \frac{\d z}{(iz)^{m+1}}
\end{align*}
to justify a deformation to a segment of the contour $\Gamma_2$. 
Call this deformed contours $\Sigma_2$, again see Figure~\ref{Figure:Sigma}.  
Cauchy's Theorem justifies adding additional contour integrals, in the upper-half plane,  which lie in $S$ and, say, tend to a steepest descent direction of $e^{-i\omega_n \sigma^n k^n}$, see Figure~\ref{Figure:AddContours}.  
The final step is to show that these additional contours can be joined with the original contour and deformed to $\cup_j \Gamma_j$.  
After some thought, it can be seen that is suffices to show that $\Gamma_1$ and $\Gamma_2$ can be deformed so that they connect with an added contour.  The results of Lemmas~\ref{Lemma:Gammaj} and \ref{Lemma:WhichGammaj} demonstrate this.  
See Figures~\ref{Figure:Join} and \ref{Figure:Gammaj} for a demonstration the deformation process.
\end{proof}
\end{lemma}

We now consider the truncation of the integration domain.
\begin{lemma}\label{Lemma:trunc}
\begin{align*}
\int_{\Gamma_j} e^{X(iz-i\omega_n \sigma^n z^n-i R_{|x|/t}(z))} \frac{\d z}{(iz)^{m+1}}  = \int_{\Gamma_j \cap \{|z| < R\}} e^{X(iz-i\omega_n \sigma^n z^n-i R_{|x|/t}(z))} &\frac{\d z}{(iz)^{m+1}} + \bigo (e^{- c_R X}),
\end{align*}
where $c_R>0$ is independent of $|x|/t$.
\begin{proof}
Consider 
\begin{align*}
I_R = \int_{\Gamma_j \cap \{|z| \geq R\}} e^{X(iz-i\omega_n \sigma^n z^n-i R_{|x|/t}(z))} \frac{\d z}{(iz)^{m+1}},
\end{align*}
and let $z = r e^{i\theta(r)}$ on one of the components of $\Gamma_j \cap \{|z| \geq R\}$. Mirroring the calculations above $\theta(r) = \theta_0 + \bigo(r^{-1})$ where $\theta_0$ is a steepest descent direction for $e^{-i\omega_n \sigma^n k^n}$ and $\theta'(r) = \bigo(r^{-1})$. Thus, $|\d z| \leq C \d r$ and
\begin{align*}
B_R = \int_{R}^\infty \left|e^{X(-i\omega_n \sigma^n r^n (\cos(n \theta(r)) +\bigo(r^{-1}))} \right|\frac{\d r}{|r|^{m+1}} = \int_{R}^\infty e^{-X|\omega_n \sigma^n| r^n (1+ |\bigo(r^{-1})|)} \frac{\d r}{|r|^{m+1}}.
\end{align*}
Let $R$ be large enough so that the $\bigo(r^{-1})$ term is less than $1/2$ (uniformly for $|x/t|$ sufficiently large).  Then
\begin{align*}
B_R \leq e^{-X |\omega_n \sigma^n| R^n/4} \int_{R}^\infty e^{-X|\omega_n \sigma^n| r^n/4} \frac{\d r}{|r|^{m+1}} \leq C_{m,R} e^{-X |\omega_n \sigma^n| R^n/4}.
\end{align*}
From this we can conclude that for fixed $m$ and fixed $R$, sufficiently large, $I_R = \bigo(e^{-X |\omega_n \sigma^n| R^n/4})$.
\end{proof}
\end{lemma}

We are now prepared to complete the steepest descent analysis.
Recall that 
\begin{align*}
I_{\omega,m}(x,t) = \frac{1}{2\pi} \sum_{j=1}^{N(n)}\int_{\Gamma_j} \frac{ e^{ikx-i\omega(k)t}}{(ik)^{m+1}} \,\d k.
\end{align*}



\begin{proof}[Proof of Theorem~\ref{Thm:Kernel}]
We perform the steepest descent analysis of
\begin{align*}
L_{m,|x|/t}(X) = \int_{\Gamma_j \cap \{|z| < R\}} e^{X(iz-i\omega_n \sigma^n z^n-i R_{|x|/t}(z))} \frac{\d z}{(iz)^{m+1}},
\end{align*}
as $X$ becomes large with $|x|/t \rightarrow \infty$.  We assume $R$ is chosen sufficiently large in the sense of Lemma~\ref{Lemma:trunc}.  We use the change of variables $k = \tau_{x/t}(s) = z_j + s v_{|x|/t}(s)$ to write
\begin{align*}
  L_{m,|x|/t}(X) &= \int_{c(R)_-}^{c(R)_+} e^{-Xs^2} F_{m,|x|/t}(s) \d s \\
  &= F_{m,|x|/t}(0)\int_{c(R)_-}^{c(R)_+} e^{-Xs^2}  \d s + F'_{m,|x|/t}(0) \int_{c(R)_-}^{c(R)_+} e^{-Xs^2} s \d s  + \half \int_{c(R)_-}^{c(R)_+} e^{-Xs^2} s^2  F''_{m,|x|/t}(\xi(s))\d s .
\end{align*}
Here $c(R)_\pm$ are chosen so that $|\tau_{|x|/t}(c(R)_\pm)| = R$.  While these functions of $R$ also depend on $x$ and $t$, it is inconsequential because they tend to finite limit as $|x|/t \rightarrow \infty$.  From Lemma~\ref{Lemma:Constants} the second integral is bounded by (assuming $c(R)_+ \geq c(R)_-$)
\begin{align*}
C  \int_{|c(R)_-|}^{c(R)_+} e^{-Xs^2} s \d s = \bigo(e^{-X |c(R)_-|}).
\end{align*}
Thus the error term is given by the third integral which by Lemma~\ref{Lemma:Constants} is $\bigo(X^{-3/2})$.  We find
\begin{align*}
L_{m,|x|/t}(X) = \sqrt{\pi} F_{m,|x|/t}(0) X^{-1/2} + \bigo(X^{-3/2}).
\end{align*}
We can confirm that the first term is indeed of higher-order because $F_{m,|x|/t}(0)$ as a definite limit as $|x|/t \rightarrow \infty$.  Computing $F_{m,|x|/t}(0)$ explicitly and using \eqref{e:rescale} we find the result.
\end{proof}

\begin{figure}[tbp]
\centering
\subfigure[]{\includegraphics[width=.49\linewidth]{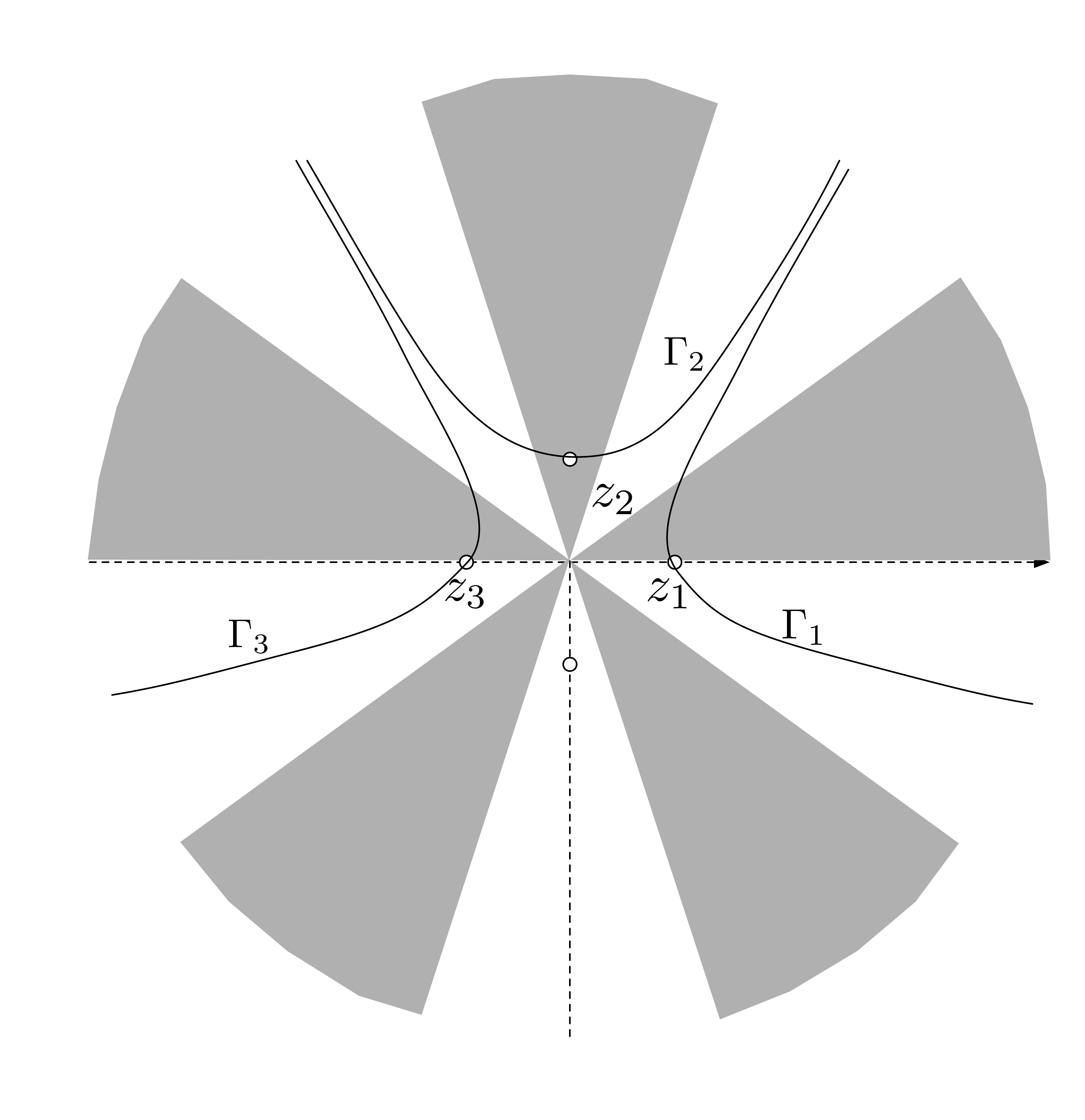}\label{Figure:GammajOdd}}
\subfigure[]{\includegraphics[width=.49\linewidth]{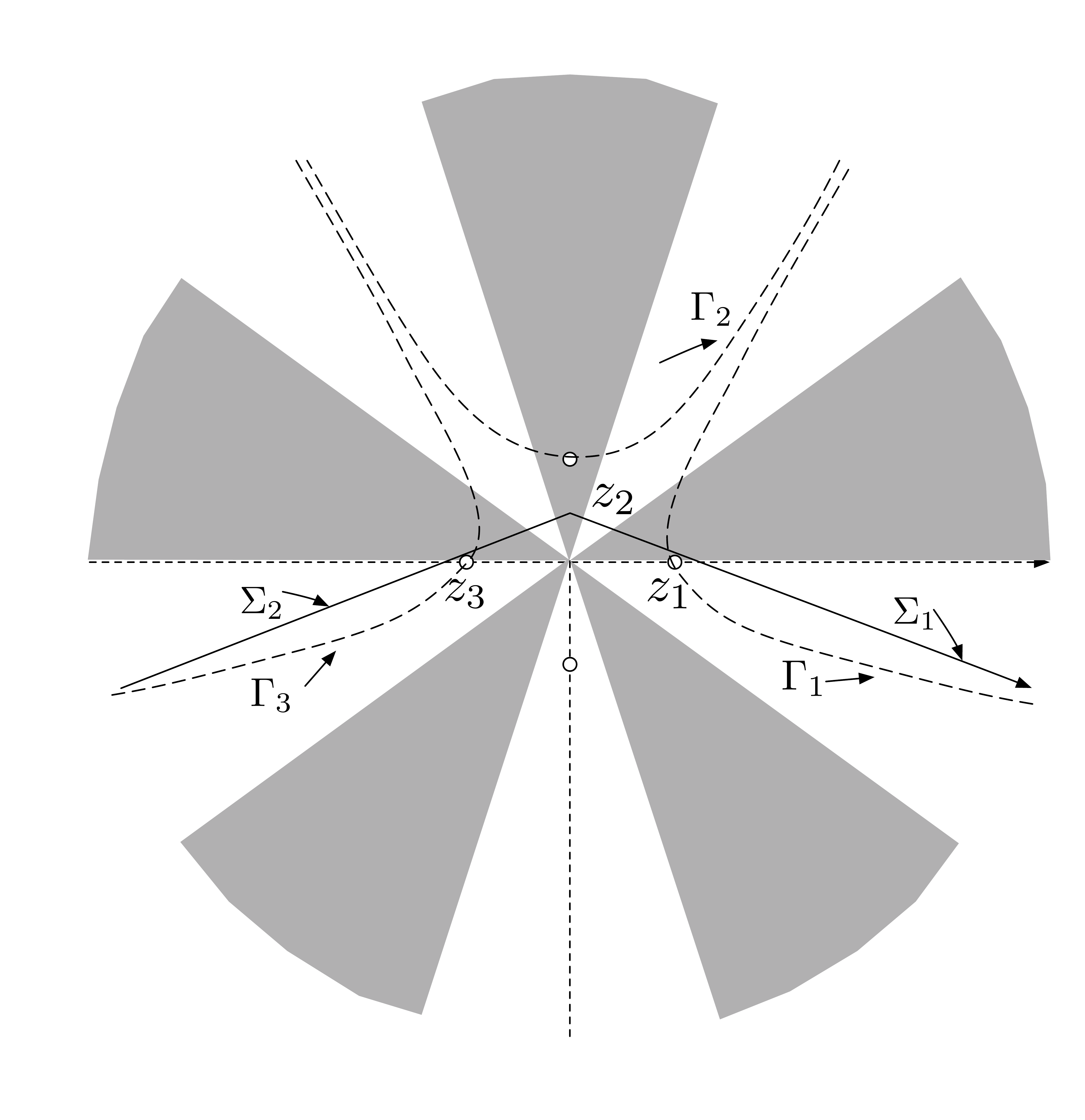}\label{Figure:SigmaOdd}}
\caption{A schematic for $\omega(k) = k^5$. The shaded region $S^c = \{z: \Im \omega_n \sigma^n z^n > 0 \}$ where the function $e^{-i\omega_n \sigma^n z^n}$ has growth. The circles represent the stationary points where $\Phi_{|x|/t}'(z)  = 0$.  (a) The contours $\Gamma_j$ which are along the global paths of steepest descent.  (b)  The initial deformation of the integral representation of $K_t(x)$ to the contours $\Sigma_1$ and $\Sigma_2$. \label{Figure:Odd1}}
\end{figure}

\begin{figure}[tbp]
\centering
\subfigure[]{\includegraphics[width=.49\linewidth]{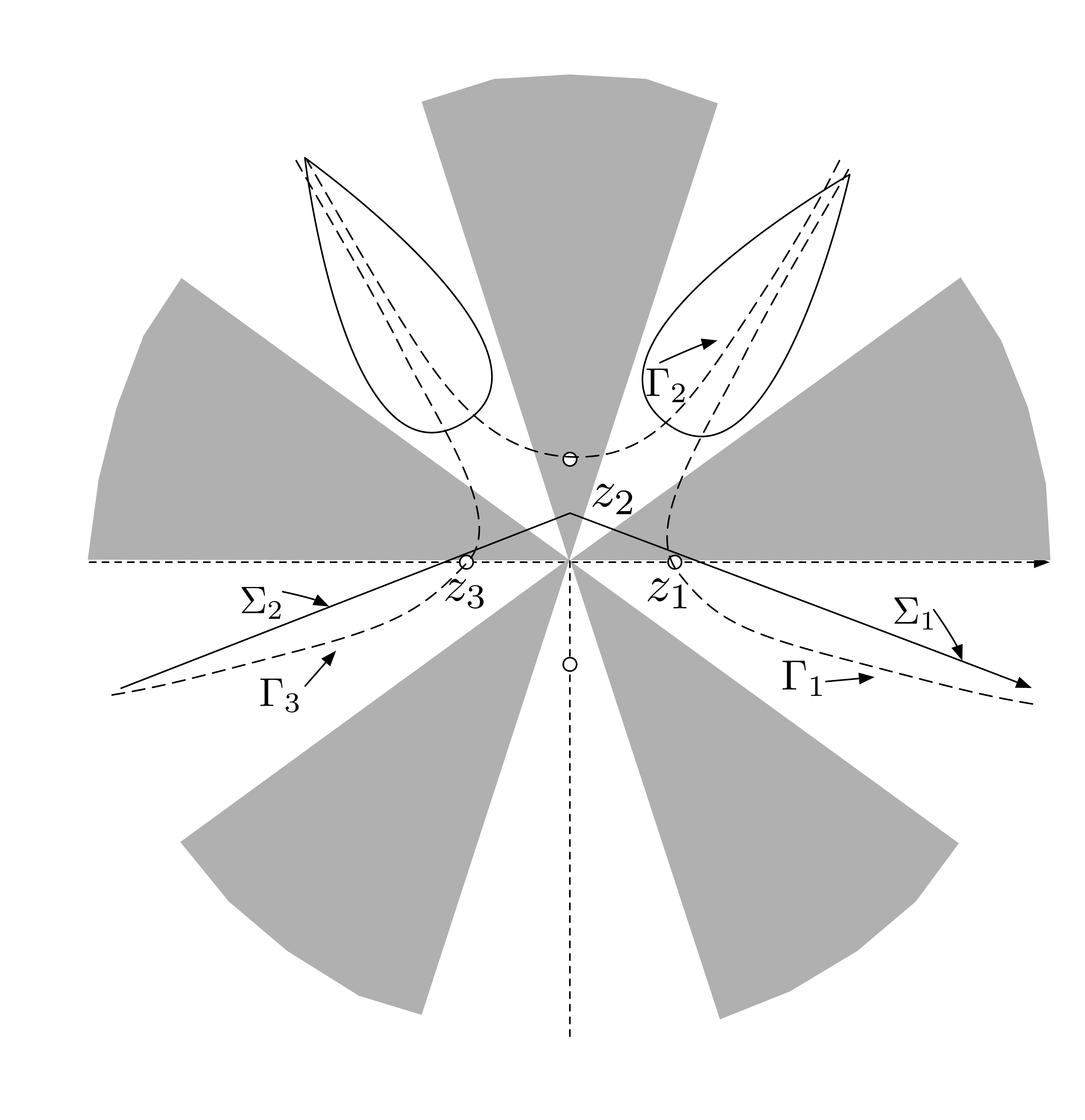}\label{Figure:AddOdd}}
\subfigure[]{\includegraphics[width=.49\linewidth]{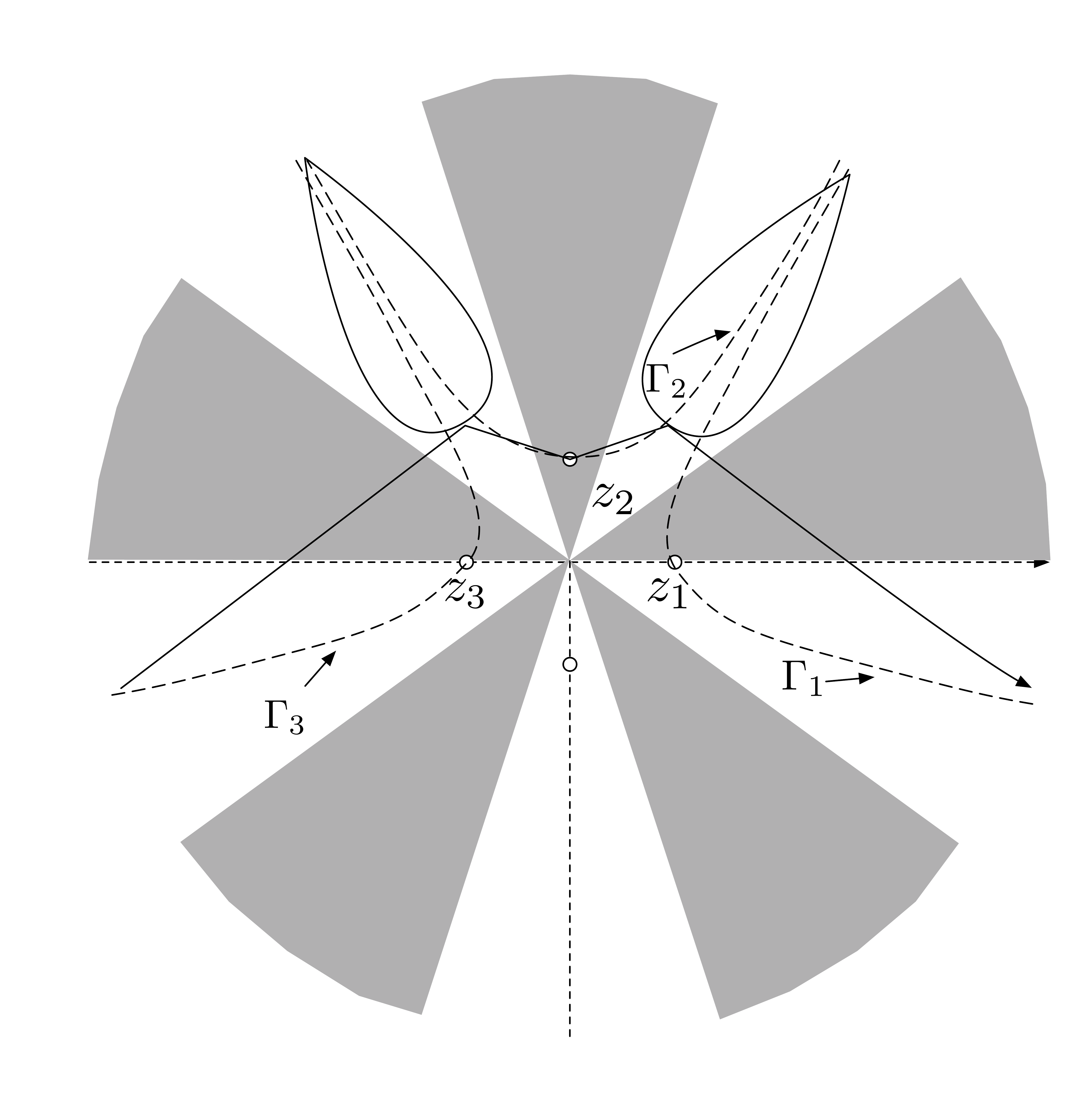}\label{Figure:JoinOdd}}
\caption{A a continuation schematic for $\omega(k) = k^5$.   The shaded region $S^c = \{z: \Im \omega_n \sigma^n z^n > 0 \}$ where the function $e^{-i\omega_n \sigma^n z^n}$ has growth. The circles represent the stationary points where $\Phi_{|x|/t}'(z) = 0$.  (a) The addition of loop contours that contribute nothing (by Cauchy's Theorem) to the integral representation for $K_t(x)$.  (b)  The deformation of $\Sigma_1$ and $\Sigma_2$ to join with the loop contours.  At this point it is clear that these contours may be deformed to $\Gamma_1 \cup \Gamma_2 \cup \Gamma_3$.\label{Figure:Odd2}}
\end{figure}

To clarify the various cases, in Figures~\ref{Figure:Odd1} and \ref{Figure:Odd2} we show a schematic for $\omega(k) = k^5$ for $x > 0$.
For $x < 0$, it suffices to consider $\omega(k) = -k^5$ with $x > 0$. 
Then all stationary points have non-zero imaginary parts and the shaded regions are switched with the unshaded regions.

We now are ready to prove the regularity theorem for linear dispersive equations, namely Theorem~\ref{t:differentiable}.

\smallskip
\begin{proof}[Proof of Theorem~\ref{t:differentiable}]~
We write $q(x,t) = K_t * q_o(x)$.
Because $q_0 \in L^2(\mathbb R)$, the weak solution is given by this convolution.  We must consider the convolution integral
\begin{align}\label{e:convo}
q(x,t) = \int_{-\infty}^{\infty} K_t(x-y)q_0(y) dy.
\end{align}
Let $B_x$ be a bounded open interval containing $x$.  
By Theorem~\ref{Thm:Kernel}
\begin{align*}
|K_t^{(m)}(x-y)q_0(y)| \leq C\left(1+|x-y|^{\frac{2m-n+2}{2(n-1)}}\right) (1+|y|)^{-\ell} (1+|y|)^\ell |q_0(y)|.
\end{align*}
For $\ell \geq \frac{2m -n + 2}{2(n-1)}$
\begin{align*}
\sup_{(x,y) \in B_x \times \mathbb R} C\left(1+|x-y|^{\frac{2m-n+2}{2(n-1)}}\right) (1+|y|)^{-\ell} = C_{\ell,B_x} < \infty.
\end{align*}
This is sufficient to justify the differentiation under the integral (see \cite[Theorem 2.27]{Folland}).  To see differentiability in $t$, we note that $\partial_t^j K_t(x)$ satisfies the same bounds, up to a constant, as $K_t^{(jn)}(x)$.
\end{proof}

Recall that a direct consequence of Theorem~\ref{t:differentiable} is Corollary~\ref{c:classical} about the requirements on the IC
to obtain a classical solution of the IVP.


%

\section{Analysis of the residual}\label{a:residual}

We now obtain estimates for the residual.  Recall the decomposition~\eref{e:qdecomp}, which we rewrite here for convenience:
\[
q(x,t) = q_c + [q_o(c)]I_{0,\omega}(x-c,t) + q_\res(x-c,t)\,, 
\]
where
\[
q_\res(y,t) = \frac1{2\pi}\int_\Real \e^{ikc}\frac{\e^{i\theta(y,t,k)}-1}{ik}\,F(k)\,\d k\,,
\]
and $\theta(x,t,k)$ was defined in~\eref{e:thetadef}.
It is trivial to see that $q_\res(0,0) = 0$.

\begin{remark}
One must keep in mind that generically, $F\not\in L^1(\mathbb R)$.  
For if it was, then $\^q_o'(k)$, defined on each interval of differentiability of $\^q_o$, would be continuous,
We did not require continuity in Assumption~\ref{Assume:One}, 
so we have no reason to assume that such a condition would be satisfied.  
This fact complicates some of the estimates that follow.
\end{remark}

We now want to show that $q_\res(x-c,t)$ is continuous as $(x,t)\to(c,0)$. 
We apply general arguments to analyze the behavior of as $(x,t) \rightarrow (c,0)$ of
\begin{gather*}
\int_{\mathbb R} e^{ikc} S_j(x-c,t,k) F(k) \, \d k, 
\\
\noalign{\noindent where}
S_j(y,t,k) = \frac1{(ik)^{j+1}}\bigg({e^{i\theta(y,t,k)} - \sum_{r = 0}^j a_r(y,t) k^r}\bigg),
\end{gather*}
and where the $a_r(y,t)$ are the Taylor coefficients for $e^{i\theta(y,t,k)}$ at $k = 0$.   
To understand the behavior of these coefficients, we write [recalling~\eref{e:dispersion}]
\begin{gather*}
\theta(y,t,k) = ky + \sum_{j=2}^{n} \omega_j (k\,t^{1/j})^j\,.
\end{gather*}
We also note that $a_r(y,t)$ is expressible as a sum of terms of the form $y^{\beta_0} \prod\nolimits_{i=1}^n t^{\beta_i/j}$, 
with $\sum\nolimits_{i}\beta_i = r$.  
Taking $|x-c|^n \leq C t$, then
\begin{gather*}
(x-c)^{\beta_0} \prod_{i=1}^n t^{\beta_i/j} = \bigo(t^{r/n}).
\end{gather*}
Thus we have $a_r(x-c,t)= \bigo(t^{r/n})$ as $t\to0$ with $|x-c|^n \leq C t$.

We now use these results to estimate the $L^q(\mathbb R)$ norm of $S_j(y,t,k)$. 
By Taylor's theorem, $S_j(y,t,k)$, on $[0,\alpha]$, 
is bounded by a polynomial of order $(n-2)(j+1)$ in $\alpha$ with coefficients that are of order $t^{(n-1)(j+1)/n}$. 
Thus for $\alpha > 1$
\begin{gather*}
\left(\int_0^\alpha |S_j(k;x,t)|^q dk\right)^{1/q} \leq C t^{(n-1)(j+1)/n} \alpha^{(n-2)(j+1)+1/q}.
\end{gather*}
Furthermore
\begin{gather*}
\left(\int_\alpha^\infty |S_j(k;x,t)|^q dk\right)^{1/q} \leq 2\sum_{r=0}^j a_r(x-c,t) \alpha^{r-j-1 + 1/q}.
\end{gather*}
We note that when $\alpha = t^{-1/n}$ both integrals are of order $t^{(j+1)/n - 1/(qn)}$. Therefore we obtain

\begin{lemma}
\label{Lemma:ResEst-HO}
Suppose $F \in L^p(\mathbb R)$ and $|x-c| \leq C|t|^n$ for $c \in \mathbb R$.  Then
\begin{align*}
\left| \int_{\mathbb R} e^{ikc} S_j(x-c,t,k) F(k) \d k \right| \leq C_{p,j,\omega} t^{j/n+1/(np)} \|F\|_{L^p(\mathbb R)}.
\end{align*}
\end{lemma}

\begin{remark}
This lemma also allows the justification of an expansion of
\begin{align*}
\int_{\mathbb R} e^{i\theta(x,t,k)} F(k) \,\d k
\end{align*}
about the point $(x,t) = (c,0)$ when $F(\cdot)(1 + |\cdot|)^{j+1} \in L^2(\mathbb R)$.  
Indeed,
\begin{align*}
\left | \int_{\mathbb R} e^{i\theta(x,t,k)} F(k) \,\d k
 - \sum_{r=0}^{j}\int_{\mathbb R} \frac{ a_r(x-c,t)k^r}{(ik)^{j+1}} (ik)^{j+1} F(k) \,\d k \right| 
 &= \left| \int_{\mathbb R} e^{ikc} S_{j}(k;x,t) (ik)^j F(k) \,\d k \right| \\
 &\leq C_{j,\omega}t^{j/n + 1/(2n)} \|F(\cdot)(1 + |\cdot|)^{j+1}\|_{L^2(\mathbb R)}.
\end{align*}
\end{remark}

\section{Approximation of ICs by smooth data}\label{a:approx}

In this appendix we discuss the behavior of the solution of \eqref{e:pde} with discontinuous data 
when it can be approximated, in an $L^1(\mathbb R)$ sense, by smooth data.  

Recall that, when the IC is continuous, the solution converges uniformly to it in the limit $t\downarrow0$,
and therefore it will not exhibit the Gibbs phenomenon as $t\to0$. 
Nonetheless, we next show that, if the IC is a small perturbation of a discontinuous function,
the solution exhibits Gibbs-like behavior \textit{at finite times}.
To see this, consider again the expression \eqref{e:convo}
\begin{align*}
q(x,t) = \int_{-\infty}^{\infty} K_t(x-y)q_0(y) dy.
\end{align*}
From \eqref{e:largespace} we know that, for $t > 0$, there exists $C_{m,\omega,t} > 0$ such that
\begin{align*}
|K_t^{(m)}(x)| \leq C_{m,\omega,t} (1+|x|)^{\frac{2m-n+2}{2(n-1)}}.
\end{align*}
Also, from Young's inequality we have for $|x|\leq R$, $R > 0$,
\begin{align*}
  |\partial_x^m q(x,t)| &\leq  C_{m,\omega,t,R} \|q_o\|_{L^1_{m,n}(\mathbb R)},\\
  \|q_o\|_{L^1_{m,n}(\mathbb R)} &\triangleq \int_{\mathbb R} |q_o(x)|(1+|x|)^{\frac{2m-n+2}{2(n-1)}} dx,
\end{align*}
with a new constant $C_{m,\omega,t,R} > 0$.  Now suppose one has a sequence $\{q_{o,\delta}\}_{\delta > 0}$ of continuous ICs which 
converges to a discontinuous function $q_o$ in the $L^1_{m,n}(\mathbb R)$ norm as $\delta \downarrow 0$. 
Let $q_\delta(x,t)$ and $q(x,t)$ be the solution of \eqref{e:pde} with initial data $q_{o,\delta}$ and $q_o$, respectively.  
It is straightforward to see that, for all fixed $t>0$ and for all $j = 0,1,\ldots,m$,
\begin{align}
|\partial_x^j(q(x,t) - q_\delta(x,t))| \goto 0\,,\qquad \delta\downarrow 0\,.
\label{e:approximation}
\end{align}
Equation~\eref{e:approximation} means that 
\textit{a Gibbs-like phenomenon similar to the one arising for $q$ as $t \downarrow 0$ will also be observed for $q_\delta$ at finite times, 
provided $\delta$ is sufficiently small}.  
Of course this statement does not hold uniformly as $t \downarrow 0$,
because $C_{m,\omega,t} = \bigo(t^{({-m-1/2})/({n-1})})$ as $t\downarrow0$ (see \eqref{e:smalltime}).

To illustrate these results, consider the following example,
in which $q_o$ is discontinuous at $x=\pm1$ but the discontinuity at $x=-1$ is smoothed out in $q_{o,\delta}$:
\begin{align}
\label{e:A4ICs}
  q_o(x) = \begin{cases} 
    1, &|x| < 1,\\
    0, \otherwise,
    \end{cases}
\quad
q_{o,\delta}(x) = \begin{cases}
    (x +1 + \delta)/\delta, &-1-\delta < x \le 1,\\
    1, &|x| < 1,\\
    0, \otherwise.
    \end{cases}
\end{align}
\begin{figure}[b!]
\centering
\includegraphics[width=\linewidth]{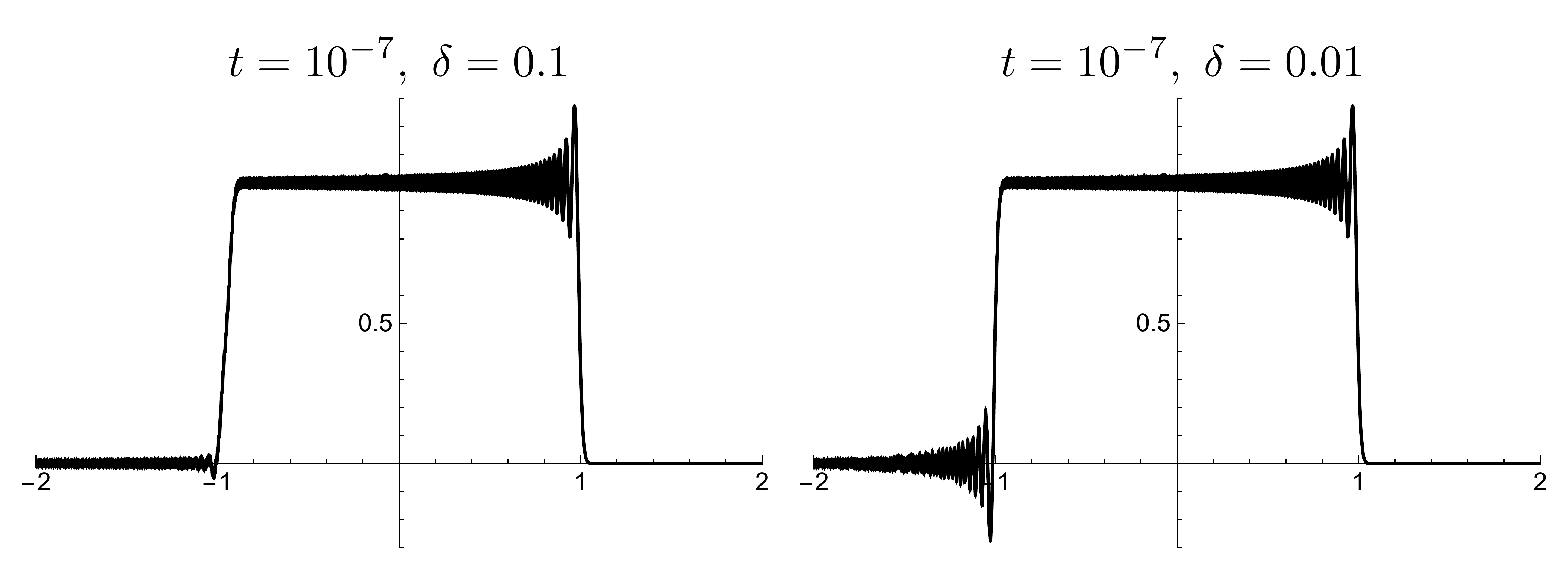}
\caption{\label{Figure:delta}
The solution $q_{\delta}(x,t)$ with $\omega(k) = -k^3$ and IC~\eref{e:A4ICs}
at $t=10^{-7}$ for two different values of $\delta$.  
Note how for $\delta = 0.1$ the Gibbs oscillations near $x=-1$ are absent, 
but $\delta = 0.01$ is sufficient for the solution to exhibit a Gibbs-like finite overshoot
even at such extremely small values of time.}
\end{figure}%
\begin{figure}[t!]
\centering
\includegraphics[width=\linewidth]{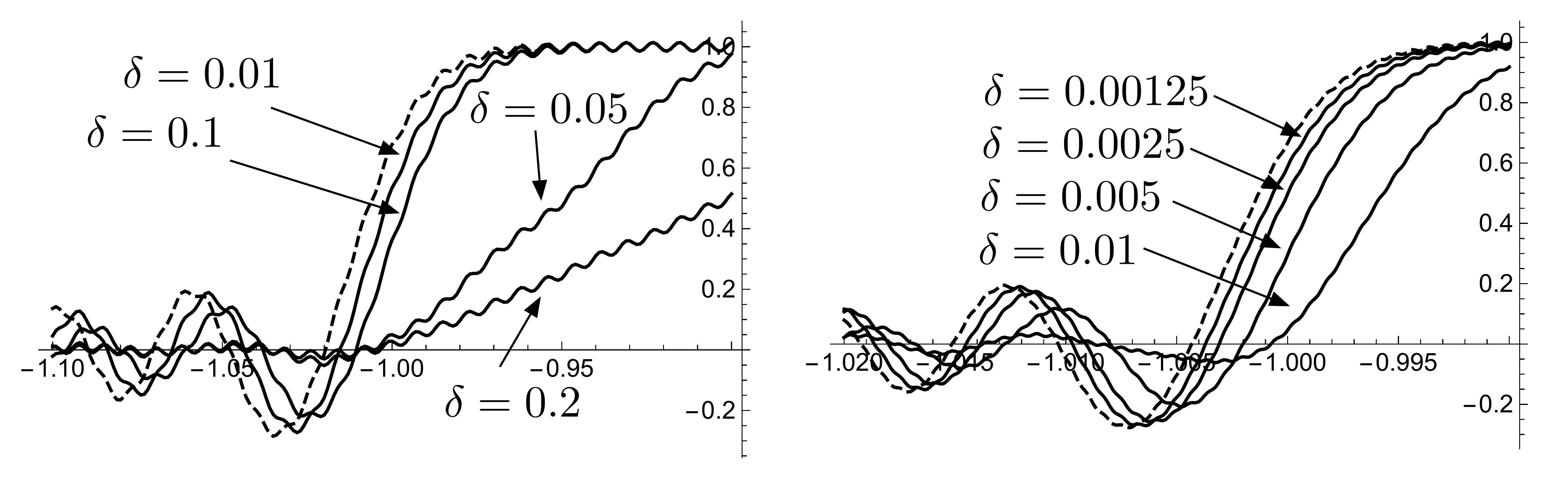}
\caption{\label{Figure:delta-zoom}
The solutions $q(x,t)$ (dashed lines) and $q_{\delta}(x,t)$ (solid lines) 
with $\omega(k) = -k^3$ and ICs~\eref{e:A4ICs} at $t = 10^{-7}$ (left) and $t=10^{-9}$ (right).  
Note how moderate values of $\delta$ make $q_\delta(x,t)$ a good approximation of $q(x,t)$ at small times.
}
\end{figure}%
Obviously $q_{o,\delta} \goto q_o$ in any $L^1_{m,n}(\mathbb R)$ norm.  
Also, the specific form of $q_o$ and $q_{o,\delta}$ ensures that the corresponding solutions
$q(x,t)$ and $q_\delta(x,t)$ are expressible in terms of the special functions $I_{\omega,m}$,
and therefore can be computed to arbitrary precision using the methods discussed earlier.
The solution behavior is displayed in Figs.~\ref{Figure:delta} and~\ref{Figure:delta-zoom}.
While $q_{\delta}(x,t)$ converges uniformly to $q_{o,\delta}(x)$ near $x = -1$ as $t \downarrow 0$, 
$q(x,t)$ does not converge uniformly to $q_{o}(x)$ near $x = -1$ as $t \downarrow 0$.
Nonetheless, for fixed $ t > 0$,  $q_{\delta}(x,t)$ converges uniformly to $q(x,t)$ near $x=-1$ as $\delta\to0$.


\section*{Acknowledgements}  
We thank Mark Ablowitz, Bernard Deconinck, Pierre Germain and Nick Trefethen for interesting discussions related to this work.  This work was partially supported by the National Science Foundation under grant numbers DMS-1311847 and DMS-1303018.


\bibliographystyle{plain}
\bigskip
\begingroup
\makeatletter
\def\@biblabel#1{#1.}
\def\journal#1{\textit{\frenchspacing #1}}
\def\title#1{``#1''}
\def\booktitle#1{\textsl{#1}}
\def\v#1{\textbf{#1}}

\endgroup

\end{document}